# High Order Numerical Scheme for Generalized Fractional Diffusion Equations

Kamlesh Kumar and Rajesh K. Pandey*[1]

*Department of Mathematics, Manav Rachna University, Haryana India

*Department of Mathematical Sciences, Indian Institute of Technology (BHU) Varanasi, Uttar Pradesh, India

## Abstract:

In this paper, a higher order finite difference scheme is proposed for Generalized Fractional Diffusion Equations (GFDEs). The fractional diffusion equation is considered in terms of the generalized fractional derivatives (GFDs) which uses the scale and weight functions in the definition. The GFD reduces to the Riemann-Liouville, Caputo derivatives and other fractional derivatives in a particular case. Due to importance of the scale and the weight functions in describing behaviour of real-life physical systems, we present the solutions of the GFDEs by considering various scale and weight functions. The convergence and stability analysis are also discussed for finite difference scheme (FDS) to validate the proposed method. We consider test examples for numerical simulation of FDS to justify the proposed numerical method.

## 1. Introduction

Fractional derivative is the generalization of integer order derivative, where order of derivative is a fraction. Due to its wide applications in solving real world problems, it has become a canter of attention for researchers working in the field of fractional partial differential equations. The subject that deals with fractional diffusion equations (FDEs) is given by replacing the integer time and/or space derivatives with fractional derivatives. Applications of FDEs include problems based on particle tracking [1], ion channel gating [2], population dynamics [3], option prices in markets [4] and other science and engineering problems namely electromagnetic [5], viscosity [6] etc. In literature, different methods for solving FDEs and other fractional models are studied since the introduction of the fractional derivatives [7, 8, 9]. Finite difference method (FDM) is one of the popular methods to provide the numerical solution with higher-order accuracy. Nowadays, initial and boundary value problems for FDEs have been discussed thoroughly in various diffusion systems [10-28].

In [29], Agrawal presented the generalized fractional derivatives (GFD) which depends upon scale and weight functions. For the special choices of scale and weight functions, GFD reduces to Caputo, Riemann–Liouville, Riesz, Hadamard and other types of fractional derivatives. The scale function is useful for adjusting simulation time via stretching and compressing for the better demonstration of solution. The weight functions allow us to assess simulation phenomena differently at different times as it extends the kernel of fractional operator to provide high level of flexibility in formulation of the problems.

[1] Corresponding Author: kkp.iitbhu@gmail.com (K. Kumar); rkpandey.mat@iitbhu.ac.in (R. K. Pandey)

We consider the time fractional generalized diffusion equation (TFGDEs) defined in terms of GFD as follows:

$$\begin{cases} {}^*\partial_t^\alpha u(x\,,t) = \delta \frac{\partial^2 u(x\,,t)}{\partial^2 x} + f(x\,,t) \quad (x\,,t) \in [a\,,\ b] \times (0\,,T], \\ u(x\,,0) = \varphi(x), \qquad x \in [a\,,\ b], \\ u(a\,,t) = f_1(t),\ u(b\,,t) = f_2(t),\ \ t > 0, \end{cases} \tag{1}$$

where $T > 0,\ \ 0 < \alpha < 1,\ \ \delta > 0$ is real numbers, $f(x\,,t)$, $f_1(t)$ and $f_2(t)$ are known functions. And, ${}^*\partial_t^\alpha u(x\,,t)$ denotes the generalized fractional derivative operator, which is defined for function $u$ with the help of scale function $z(t)$ and weight function $w(t)$ as,

$${}^*\partial_t^\alpha u(x\,,t) = \frac{[w(t)]^{-1}}{\Gamma(1-\alpha)} \int_0^t \frac{\left(w(\tau)u(x\,,\tau)\right)'}{\left(z(t)-z(\tau)\right)^\alpha} d\tau, \qquad 0 < \alpha < 1. \tag{2}$$

Due to presence of scale and weight functions, it is difficult to find analytical solution of FDEs. Hence numerical solution is necessary in order to solve such type of diffusion equations. In literature, very few research works are available related with numerical solution of GFDEs [30-35] with higher order of convergence. Motivated by the research works presented in [30-34], we present the numerical solution of TFGDEs as given in Eq. (1) with high order of convergence.

Rest of the work is structured in following manner: In Section 2, problem formulation as well as computational algorithm is presented for TFGDEs using finite difference scheme. Section 3 depicts the approximation of the GFD by using interpolating polynomials and order of error of approximation. The stability analysis and order of convergence of proposed FDS are derived in Section 4. In Section 5, we present some illustrative examples to validate the proposed numerical scheme for TFGDEs by considering different scale and weight functions. Section 6 concludes our work.

## 2. Computational Algorithm of Problem

To find the numerical solution firstly, we divide space domain $[a,\ b]$ and time domain $[0,\ T]$ into $N$ and $M$ equal parts, respectively. Suppose $\mathcal{P}_N = \{a = x_0 < x_1 < x_2 < \cdots < x_N = b\}$ be the partition of $[a,\ b]$ into $N$ equal parts with spacing $\Delta x = \frac{b-a}{N}$, and $\mathcal{P}_M = \{0 = t_0 < t_1 < t_2 < \cdots < t_M = T\}$ be the partition of $[0,\ T]$ into $M$ equal parts with spacing $\Delta t = \frac{T}{M}$. For convenience, we denote $u(x_i\,, t_j) = u_j^i$, $w(t_j) = w_j$ and $z(t_j) = z_j$. We assume $w(t) > 0$, and $z(t)$ is strictly increasing then $s = z^{-1}(v)$ by taking transformation $v = z(s)$. For numerical solution, we approximate each term of GFDE given by Eq. (1). The Caputo type generalized fractional derivative of function $u(x\,,t)$ can be approximated at node point $(x_i\,,t_{j+1})$. Firstly, we approximate of first term of Eq. (1) as follow,

$$[{}^*\partial_t^\alpha u(x\,,t)]_{(x_i,\;t_{j+1})} = \frac{[w(t_{j+1})]^{-1}}{\Gamma(1-\alpha)}\int_0^{t_{j+1}} \frac{\frac{\partial}{\partial s}[w(s)u(x_i\,,s)]}{[z(t_{j+1})-z(s)]^\alpha}ds = \frac{[w(t_{j+1})]^{-1}}{\Gamma(1-\alpha)}\sum_{k=0}^{j}\int_{t_k}^{t_{k+1}} \frac{\frac{\partial}{\partial s}[w(s)u(x_i\,,s)]}{[z(t_{j+1})-z(s)]^\alpha}ds$$

$$= \frac{[w(t_{j+1})]^{-1}}{\Gamma(1-\alpha)}\int_{t_0}^{t_1} \frac{\frac{\partial}{\partial s}[w(s)u(x_i\,,s)]}{[z(t_{j+1})-z(s)]^\alpha}ds + \frac{[w(t_{j+1})]^{-1}}{\Gamma(1-\alpha)}\sum_{k=1}^{j}\int_{t_k}^{t_{k+1}} \frac{\frac{\partial}{\partial s}[w(s)u(x_i\,,s)]}{[z(t_{j+1})-z(s)]^\alpha}ds \tag{3}$$

$$= \frac{[w(t_{j+1})]^{-1}}{\Gamma(1-\alpha)}\int_{z_0}^{z_1} \frac{1}{[z(t_{j+1})-v]^\alpha}\frac{d\left[w\left(z^{-1}(v)\right)u\left(x_i\,,z^{-1}(v)\right)\right]}{dz^{-1}(v)}dz^{-1}(v) +$$

$$\frac{[w(t_{j+1})]^{-1}}{\Gamma(1-\alpha)}\sum_{k=1}^{j}\int_{z_k}^{z_{k+1}} \frac{1}{[z(t_{j+1})-v]^\alpha}\frac{d\left[w\left(z^{-1}(v)\right)u\left(x_i\,,z^{-1}(v)\right)\right]}{dz^{-1}(v)}dz^{-1}(v) \tag{4}$$

On the first interval $[z_0\;,z_1]$, the linear interpolation $\left(p_u^1(v)\right)'$ is used for approximating the unknown function and for the other subintervals $(k \geq 1)$, the quadratic interpolation function $\left(p_u^2(v)\right)'$ for the three points $\left(z_{k-1}\,,u_{k-1}^i\right)$, $\left(z_k\,,u_k^i\right)$, $\left(z_{k+1}\,,u_{k+1}^i\right)$ is applied such that,

$$\left(p_u^1(v)\right)' = \frac{w_1u_1^i - w_0u_0^i}{z_1 - z_0}\,, \tag{5}$$

$$\left(p_u^2(v)\right)'$$
$$= \left(\frac{(v-z_k)(v-z_{k+1})}{(z_k-z_{k-1})(z_{k+1}-z_{k-1})}w_{k-1}u_{k-1}^i - \frac{(v-z_{k-1})(v-z_{k+1})}{(z_k-z_{k-1})(z_{k+1}-z_k)}w_ku_k^i + \frac{(v-z_{k-1})(v-z_k)}{(z_{k+1}-z_{k-1})(z_{k+1}-z_k)}w_{k+1}u_{k+1}^i\right)'. \tag{6}$$

$$\left(p_u^2(v)\right)'$$
$$= \frac{(2v-z_k-z_{k+1})}{(z_k-z_{k-1})(z_{k+1}-z_{k-1})}w_{k-1}u_{k-1}^i - \frac{(2v-z_{k-1}-z_{k+1})}{(z_k-z_{k-1})(z_{k+1}-z_k)}w_ku_k^i + \frac{(2v-z_{k-1}-z_k)}{(z_{k+1}-z_{k-1})(z_{k+1}-z_k)}w_{k+1}u_{k+1}^i. \tag{7}$$

Again, we consider first term of Eq. (4)

$$\frac{[w(t_{j+1})]^{-1}}{\Gamma(1-\alpha)}\int_{z_0}^{z_1} \frac{1}{[z(t_{j+1})-v]^\alpha}\frac{d\left[w\left(z^{-1}(v)\right)u\left(x_i\,,z^{-1}(v)\right)\right]}{dz^{-1}(v)}dz^{-1}(v)$$

$$= \frac{[w(t_{j+1})]^{-1}}{\Gamma(2-\alpha)}\left[\frac{w_1u_1^i - w_0u_0^i}{z_1-z_0}\right]\left[\left(z_{j+1}-z_0\right)^{1-\alpha} - \left(z_{j+1}-z_1\right)^{1-\alpha}\right]$$

$$= \frac{[w(t_{j+1})]^{-1}}{\Gamma(2-\alpha)}\left[\frac{\left(z_{j+1}-z_0\right)^{1-\alpha} - \left(z_{j+1}-z_1\right)^{1-\alpha}}{z_1-z_0}\right]w_1u_1^i - \frac{[w(t_{j+1})]^{-1}}{\Gamma(2-\alpha)}\left[\frac{\left(z_{j+1}-z_0\right)^{1-\alpha} - \left(z_{j+1}-z_1\right)^{1-\alpha}}{z_1-z_0}\right]w_0u_0^i\,. \tag{8}$$

And second term of Eq. (4)

$$\frac{[w(t_{j+1})]^{-1}}{\Gamma(1-\alpha)}\sum_{k=1}^{j}\int_{z_k}^{z_{k+1}} \frac{1}{[z(t_{j+1})-v]^\alpha}\frac{d\left[w\left(z^{-1}(v)\right)u\left(x_i\,,z^{-1}(v)\right)\right]}{dz^{-1}(v)}dz^{-1}(v)$$

$$= \sum_{k=1}^{j} a_k^j u_{k-1}^i - b_k^j u_k^i + c_k^j u_{k+1}^i \tag{9}$$

$$a_k^j = \frac{w_{k-1}}{(z_k - z_{k-1})(z_{k+1} - z_{k-1})} \frac{[w(t_{j+1})]^{-1}}{\Gamma(1-\alpha)} \left[\frac{(2z_{j+1} - z_k - z_{k+1})}{1-\alpha} p_k^j - \frac{2}{2-\alpha} q_k^j\right] . \tag{10}$$

$$b_k^j = \frac{w_k}{(z_k - z_{k-1})(z_{k+1} - z_k)} \frac{[w(t_{j+1})]^{-1}}{\Gamma(1-\alpha)} \left[\frac{(2z_{j+1} - z_{k-1} - z_{k+1})}{1-\alpha} p_k^j - \frac{2}{2-\alpha} q_k^j\right] . \tag{11}$$

$$c_k^j = \frac{w_{k+1}}{(z_{k+1} - z_{k-1})(z_{k+1} - z_k)} \frac{[w(t_{j+1})]^{-1}}{\Gamma(1-\alpha)} \left[\frac{(2z_{j+1} - z_{k-1} - z_k)}{1-\alpha} p_k^j - \frac{2}{2-\alpha} q_k^j\right] . \tag{12}$$

$$p_k^j = \left(z_{j+1} - z_k\right)^{1-\alpha} - \left(z_{j+1} - z_{k+1}\right)^{1-\alpha} . \tag{13}$$

$$q_k^j = \left(z_{j+1} - z_k\right)^{2-\alpha} - \left(z_{j+1} - z_{k+1}\right)^{2-\alpha} . \tag{14}$$

Finally, we get approximation of first term of Eq. (1) as,

$$[{}^*\partial_t^\alpha u(x\,,t)]_{(x_i\, t_{j+1})} \approx \frac{[w(t_{j+1})]^{-1}}{\Gamma(2-\alpha)} \left[\frac{p_0^j}{z_1 - z_0}\right] \left[w_1 u_1^i - w_0 u_0^i\right] + \sum_{k=1}^{j} \left(a_k^j u_{k-1}^i - b_k^j u_k^i + c_k^j u_{k+1}^i\right)$$

$$= p_j u_1^i - q_j u_0^i + \sum_{k=1}^{j} \left(a_k^j u_{k-1}^i - b_k^j u_k^i + c_k^j u_{k+1}^i\right) , \tag{15}$$

where,

$$p_j = \frac{[w(t_{j+1})]^{-1}}{\Gamma(2-\alpha)} \left[\frac{p_0^j w_1}{z_1 - z_0}\right] \tag{16}$$

$$q_j = \frac{[w(t_{j+1})]^{-1}}{\Gamma(2-\alpha)} \left[\frac{p_0^j w_0}{z_1 - z_0}\right] \tag{17}$$

For second order derivative of $u(x\,,t)$ in the spatial direction to appear in Eq. (1), we use approximation as

$$\left[\frac{\partial^2 u(x\,,t)}{\partial^2 x}\right]_{(x_i\, t_{j+1})} \approx \frac{u(x_{i+1}\,,t_{j+1}) - 2u(x_i\,,t_{j+1}) - u(x_{i-1}\,,t_{j+1})}{(\Delta x)^2} = \frac{u_{j+1}^{i+1} - 2u_{j+1}^i + u_{j+1}^{i-1}}{(\Delta x)^2} . \tag{18}$$

Hence, from Eq. (15), (18) and Eq. (1), the scheme takes the form:

$$p_j u_1^i - q_j u_0^i + \sum_{k=1}^{j} \left(a_k^j u_{k-1}^i - b_k^j u_k^i + c_k^j u_{k+1}^i\right) = \delta \left[\frac{u_{j+1}^{i+1} - 2u_{j+1}^i + u_{j+1}^{i-1}}{(\Delta x)^2}\right] + f_{j+1}^i . \tag{19}$$

For simplicity point of view Eq. (19) can be think as,

$$\mu\left(u_{j+1}^{i+1} - 2u_{j+1}^i + u_{j+1}^{i-1}\right) = p_j u_1^i - q_j u_0^i + \sum_{k=1}^{j} \left(a_k^j u_{k-1}^i - b_k^j u_k^i + c_k^j u_{k+1}^i\right) - f_{j+1}^i \tag{20}$$

where,

$$\mu = \frac{\delta}{(\Delta x)^2} \ , \tag{21}$$

$$\mu u_{j+1}^{i+1} - \left(2\mu + c_j^j\right) u_{j+1}^i + \mu u_{j+1}^{i-1} = p_j u_1^i - q_j u_0^i + \left(a_j^j u_{j-1}^i - b_j^j u_j^i\right) + \sum_{k=1}^{j-1}\left(a_k^j u_{k-1}^i - b_k^j u_k^i + c_k^j u_{k+1}^i\right) - f_{j+1}^i \ , \tag{22}$$

where, $0 \le j \le M-1$ and $1 \le i \le N-1$.

Using, $L_j^i = -\left(2\mu + c_j^j\right)$, Eq. (19) takes the matrix form for computation of $u_{j+1}^i$ , $j = 0, 1, 2, \cdots, M-1$,

$$K_{j+1} U_{j+1} = F_{j+1} \ , \quad 0 \le j \le M-1, \tag{23}$$

where,

$$K_{j+1} = \begin{pmatrix} L_j^1 & \mu & & & & \\ \mu & L_j^2 & \mu & & & \\ & \ddots & \ddots & \ddots & & \\ & & & \mu & L_j^{N-2} & \mu \\ & & & & \mu & L_j^{N-1} \end{pmatrix}, \tag{24}$$

$$U_{j+1} = \left[u_{j+1}^1, u_{j+1}^2, \cdots \ , u_{j+1}^i, \cdots, u_{j+1}^{N-1}\ \right]^T, \tag{25}$$

$$F_{j+1} = \left[F_{j+1}^1, F_{j+1}^2, \cdots \ , F_{j+1}^i, \cdots, F_{j+1}^{N-1}\ \right]^T, \tag{26}$$

and,

$$F_{j+1}^i = \begin{cases} -q_0 u_0^i - f_1^i \ , & \text{for } j = 0 \\ p_j u_1^i - q_j u_0^i + \left(a_j^j u_{j-1}^i - b_j^j u_j^i\right) + \cdots & \\ + \sum_{k=1}^{j-1}\left(a_k^j u_{k-1}^i - b_k^j u_k^i + c_k^j u_{k+1}^i\right) - f_{j+1}^i \ , & \text{for } 1 \le j \le M-1 \end{cases} . \tag{27}$$

## 3. Order of Approximation for the Generalized Fractional Derivative

To analyze error of approximation for generalized time fractional derivative, we use notation for simplicity, let $g(s) = w(s) f(s)$ and $s = z^{-1}(v)$ . Hence, we find $g(v) = w\left(z^{-1}(v)\right) f\left(z^{-1}(v)\right)$, $w(t_k) = w_k$ and $z(t_k) = z_k$. According to Newton interpolation method, we interpolate $g(v)$ on $[z_0, z_1]$ and $[z_{k-1}, z_{k+1}]$ by $p_g^1(v)$ and $p_g^2(v)$ respectively such that

$$g(v) - p_g^1(v) = \frac{g''(\xi_1)}{2!}(v - z_0)(v - z_1), \qquad \xi_1 \in [z_0, z_1] \ . \tag{28}$$

$$g(v) - p_g^2(v) = \frac{g'''(\xi_2)}{3!}(v - z_{k-1})(v - z_k)(v - z_{k+1}), \qquad \xi_2 \in [z_{k-1}, z_{k+1}] \ , \tag{29}$$

where $p_g^1(v)$ and $p_g^2(v)$ are piecewise linear and quadratic interpolation polynomials using the nodes $(z_0\,,g_0)$ , $(z_1\,,g_1)$ and $(z_{k-1}\,,g_{k-1})$ $(z_k\,,g_k)$ $(z_{k+1}\,,g_{k+1})$ respectively.

Let $E^{j+1}$ be the truncation error given by

$$E^{j+1} = \frac{[w_{j+1}]^{-1}}{\Gamma(1-\alpha)}\int_{z_0}^{z_1}\frac{[g(v)-p_g^1(v)]'}{[z_{j+1}-v]^{\alpha}}dv + \frac{[w_{j+1}]^{-1}}{\Gamma(1-\alpha)}\sum_{k=1}^{j}\int_{z_k}^{z_{k+1}}\frac{1}{[z_{j+1}-v]^{\alpha}}\left[g(v)-p_g^2(v)\right]'dv. \tag{30}$$

Consider first term of Eq. (30)

$$\frac{[w_{j+1}]^{-1}}{\Gamma(1-\alpha)}\int_{z_0}^{z_1}\frac{1}{[z_{j+1}-v]^{\alpha}}\left[g(v)-p_g^1(v)\right]'dv = \frac{[w_{j+1}]^{-1}}{2\Gamma(1-\alpha)}\int_{z_0}^{z_1}\frac{g''(\xi_1)}{[z_{j+1}-v]^{\alpha}}[(v-z_0)(v-z_1)]'dv. \tag{31}$$

Using integration by parts, we get

$$\frac{[w_{j+1}]^{-1}}{2\Gamma(1-\alpha)}\int_{z_0}^{z_1} g''(\xi_1)\left(z_{j+1}-v\right)^{-\alpha-1}(v-z_0)(v-z_1)dv$$

$$\leq \frac{[w_{j+1}]^{-1}}{2\Gamma(1-\alpha)}\max_{t_0\leq\xi_1\leq t_{j+1}} g''(\xi_1)\left(z_{j+1}-z_1\right)^{-\alpha-1}\int_{z_0}^{z_1}(v-z_0)(z_1-v)dv$$

$$\leq \frac{-\alpha\max_{t_0\leq\xi_1\leq t_{j+1}} g''(\xi_1)\left(z_{j+1}-z_1\right)^{-\alpha-1}(z_1-z_0)^3}{w_{j+1}\Gamma(1-\alpha)\,12}. \tag{32}$$

Consider the second term of Eq. (30)

$$\frac{[w_{j+1}]^{-1}}{\Gamma(1-\alpha)}\sum_{k=1}^{j}\int_{z_k}^{z_{k+1}}\frac{1}{[z_{j+1}-v]^{\alpha}}\left[g(v)-p_g^2(v)\right]'dv.$$

Since, $g(v)-p_g^2(v) = \frac{g'''(\xi_2)}{3!}(v-z_{k-1})(v-z_k)(v-z_{k+1})$ and applying integration by parts, we get

$$\frac{\alpha[w_{j+1}]^{-1}}{\Gamma(1-\alpha)}\sum_{k=1}^{j}\int_{z_k}^{z_{k+1}}\frac{g'''(\xi_2)}{3!}(v-z_{k-1})(v-z_k)(v-z_{k+1})\left[z(t_{j+1})-v\right]^{-\alpha-1}dv$$

$$= \frac{\alpha[w_{j+1}]^{-1}}{\Gamma(1-\alpha)}\sum_{k=1}^{j-1}\int_{z_k}^{z_{k+1}}\frac{g'''(\xi_2)}{3!}(v-z_{k-1})(v-z_k)(v-z_{k+1})\left[z(t_{j+1})-v\right]^{-\alpha-1}dv$$

$$+\frac{\alpha[w_{j+1}]^{-1}}{\Gamma(1-\alpha)}\int_{z_j}^{z_{j+1}}\frac{g'''(\xi_2)}{3!}\left(v-z_{j-1}\right)\left(v-z_j\right)\left(v-z_{j+1}\right)\left[z(t_{j+1})-v\right]^{-\alpha-1}dv \tag{33}$$

Consider the first part of Eq. (33)

$$\frac{\alpha[w_{j+1}]^{-1}}{\Gamma(1-\alpha)}\sum_{k=1}^{j-1}\int_{z_k}^{z_{k+1}}\frac{g'''(\xi_2)}{3!}(v-z_{k-1})(v-z_k)(v-z_{k+1})\left[z(t_{j+1})-v\right]^{-\alpha-1}dv$$

$$\leq \frac{\alpha[w_{j+1}]^{-1}}{\Gamma(1-\alpha)}\varphi(z_{k-1}\,,z_k\,,z_{k+1})\sum_{k=1}^{j-1}\int_{z_k}^{z_{k+1}}\frac{g'''(\xi_2)}{3!}\left[z(t_{j+1})-v\right]^{-\alpha-1}dv \tag{34}$$

$$= \frac{\alpha[w_{j+1}]^{-1}}{\Gamma(1-\alpha)}\varphi(z_{k-1}\,,z_k\,,z_{k+1})\frac{\max_{t_0\leq\xi_2\leq t_j}|g'''(\xi_2)|}{3!}\int_{z_1}^{z_j}\left[z(t_{j+1})-v\right]^{-\alpha-1}dv, \tag{35}$$

where,

$$\varphi(z_{k-1}, z_k, z_{k+1}) = \frac{1}{27}\varphi_1(z_{k-1}, z_k, z_{k+1})\varphi_2(z_{k-1}, z_k, z_{k+1})\varphi_3(z_{k-1}, z_k, z_{k+1}) \tag{36}$$

$$\varphi_1(z_{k-1}, z_k, z_{k+1}) = (z_k + z_{k+1} - 2z_{k-1}) - \sigma(z_{k-1}, z_k, z_{k+1}) \tag{37}$$

$$\varphi_2(z_{k-1}, z_k, z_{k+1}) = (z_{k+1} + z_{k-1} - 2z_k) - \sigma(z_{k-1}, z_k, z_{k+1}) \tag{38}$$

$$\varphi_3(z_{k-1}, z_k, z_{k+1}) = (z_{k-1} + z_k - 2z_{k+1}) - \sigma(z_{k-1}, z_k, z_{k+1}) \tag{39}$$

$$\sigma(z_{k-1}, z_k, z_{k+1}) = \sqrt{z_{k-1}(z_{k-1} - z_k) + z_k(z_k - z_{k+1}) + z_{k+1}(z_{k+1} - z_{k-1})} \tag{40}$$

Again,

$$\int_{z_1}^{z_j}\left[z(t_{j+1}) - v\right]^{-\alpha-1} dv = \frac{1}{\alpha}\left[\left(z_{j+1} - z_j\right)^{-\alpha} - \left(z_{j+1} - z_1\right)^{-\alpha}\right] \le \left(z_{j+1} - z_j\right)^{-\alpha}. \tag{41}$$

Hence, from Eq. (35) and (41), we have first term of Eq. (33) is less than or equal to

$$\frac{\alpha[w_{j+1}]^{-1}}{\Gamma(1-\alpha)}\varphi(z_{k-1}, z_k, z_{k+1})\frac{\max_{t_0 \le \xi_2 \le t_j}|g'''(\xi_2)|}{3!}\left(z_{j+1} - z_j\right)^{-\alpha}\ . \tag{42}$$

Consider the second term of Eq. (33)

$$\frac{\alpha[w_{j+1}]^{-1}}{\Gamma(1-\alpha)}\int_{z_j}^{z_{j+1}}\frac{g'''(\xi_2)}{3!}\left(v - z_{j-1}\right)\left(v - z_j\right)\left(v - z_{j+1}\right)\left[z(t_{j+1}) - v\right]^{-\alpha-1} dv$$

$$= -\frac{\alpha[w_{j+1}]^{-1}}{3!\Gamma(1-\alpha)}\int_{z_j}^{z_{j+1}}\frac{g'''(\xi_2)}{1}\left(v - z_{j-1}\right)\left(v - z_j\right)\left[z(t_{j+1}) - v\right]^{-\alpha} dv \tag{43}$$

$$= -\frac{\alpha[w_{j+1}]^{-1}\max_{t_0 \le \xi_2 \le t_j}|g'''(\xi_2)|\left(z_{j+1}-z_j\right)^{(2-\alpha)}}{3!\Gamma(1-\alpha)}\left[\frac{(z_j - z_{j-1})}{(1-\alpha)(2-\alpha)} + \frac{2(z_{j+1}-z_j)}{(1-\alpha)(2-\alpha)(3-\alpha)}\right] \tag{44}$$

$$= -\frac{\alpha[w_{j+1}]^{-1}}{3\Gamma(1-\alpha)}\frac{\max_{t_0 \le \xi_2 \le t_j}|g'''(\xi_2)|\left(z_{j+1}-z_j\right)^{(2-\alpha)}}{(1-\alpha)(2-\alpha)}\left[\frac{(z_j - z_{j-1})}{2} + \frac{(z_{j+1}-z_j)}{(3-\alpha)}\right] \tag{45}$$

From Eqs. (32), (42) and (45), we have

$$E^{j+1} \le \frac{-\alpha \max_{t_0 \le \xi_1 \le t_{j+1}} g''(\xi_1)\left(z_{j+1}-z_1\right)^{-\alpha-1}(z_1 - z_0)^3}{w_{j+1}\Gamma(1-\alpha)\ 12}$$

$$+\frac{\alpha[w_{j+1}]^{-1}}{\Gamma(1-\alpha)}\varphi(z_{k-1}, z_k, z_{k+1})\frac{\max_{t_0 \le \xi_2 \le t_j}|g'''(\xi_2)|}{3!}\left(z_{j+1} - z_j\right)^{-\alpha}$$

$$-\frac{\alpha[w_{j+1}]^{-1}}{3\Gamma(1-\alpha)}\frac{\max_{t_0 \le \xi_2 \le t_j}|g'''(\xi_2)|\left(z_{j+1}-z_j\right)^{(2-\alpha)}}{(1-\alpha)(2-\alpha)}\left[\frac{(z_j - z_{j-1})}{2} + \frac{(z_{j+1}-z_j)}{(3-\alpha)}\right]. \tag{46}$$

**Lemma 1:** If scale function $z(t)$ is satisfying the Lipschitz condition with constant $L$ that is $(z_{j+1} - z_j) \leq L\,\Delta t$ then we have following relation

(i) $\sigma(z_{k-1}, z_k, z_{k+1}) \leq \sqrt{3}\,L\,\Delta t$

(ii) $\varphi_1(z_{k-1}, z_k, z_{k+1}) \leq L\,\Delta t(3-\sqrt{3})$

(iii) $\varphi_2(z_{k-1}, z_k, z_{k+1}) \leq -L\,\Delta t(\sqrt{3})$

(iv) $\varphi_3(z_{k-1}, z_k, z_{k+1}) \leq L\,\Delta t(-3-\sqrt{3})$

(v) $\varphi(z_{k-1}, z_k, z_{k+1}) \leq \left(\frac{2\sqrt{3}}{9}\right) L\,(\Delta t)^3 \leq \frac{1}{2} L\,(\Delta t)^3.$

It is clear from above Lemma 1 and Eq. (46), the truncation error $E^{j+1}$ achieves the order of convergence $(\Delta t)^{3-\alpha}$. If the scale function is satisfying the Lipschitz condition with constant $L$ then the truncation error $E^{j+1}$ has the form

$$E^{j+1} \leq \frac{-\alpha \max_{t_0 \leq \xi_1 \leq t_{j+1}} g''(\xi_1)(t_{j+1}-t_1)^{-\alpha-1} L^{-\alpha}(\Delta t)^3}{w_{j+1}\Gamma(1-\alpha)\,12} + \frac{\alpha[w_{j+1}]^{-1}}{\Gamma(1-\alpha)}\frac{1}{12} L^{3-\alpha} \max_{t_0 \leq \xi_2 \leq t_j} |g'''(\xi_2)|(\Delta t)^{3-\alpha}$$

$$-\frac{\alpha[w_{j+1}]^{-1}}{3\Gamma(1-\alpha)} \frac{\max_{t_0 \leq \xi_2 \leq t_j} |g'''(\xi_2)|}{(1-\alpha)(2-\alpha)} \left[\frac{1}{2} + \frac{2}{(3-\alpha)}\right] L^{3-\alpha}(\Delta t)^{3-\alpha}\,. \tag{47}$$

**Lemma 2:** Suppose $z(t)$ is strictly increasing with $z(t) \geq 0$ and $w(t)$ is increasing with $w(t) > 0$. Then the following conditions hold for $\alpha \in (0, 1)$,

(i) $A_k^j = \alpha\left[(z_{j+1} - z_k)^{2-\alpha} - (z_{j+1} - z_{k+1})^{2-\alpha}\right] - (2-\alpha)\left\{(z_{j+1} - z_k)(z_{j+1} - z_{k+1})^{1-\alpha} - (z_{j+1} - z_{k+1})(z_{j+1} - z_k)^{1-\alpha}\right\}$,

(ii) $A_k^j = 2\left[(z_{j+1} - z_k)^{2-\alpha} - (z_{j+1} - z_{k+1})^{2-\alpha}\right]$
$-(2-\alpha)(z_{k+1} - z_k)\left[(z_{j+1} - z_k)^{1-\alpha} + (z_{j+1} - z_{k+1})^{1-\alpha}\right]$,

(iii) $A_k^j > 0$,

(iv) $C_k^j > B_k^j > A_k^j > 0$,

(v) $a_k^j$, $b_k^j$ and $c_k^j$ are positive,

where

$$a_k^j = \frac{w_{k-1}}{(z_k - z_{k-1})(z_{k+1} - z_{k-1})} \frac{[w(t_{j+1})]^{-1}}{\Gamma(3-\alpha)} A_k^j \tag{48}$$

$$b_k^j = \frac{w_k}{(z_k - z_{k-1})(z_{k+1} - z_k)} \frac{[w(t_{j+1})]^{-1}}{\Gamma(1-\alpha)} B_k^j \tag{49}$$

$$c_k^j = \frac{w_{k+1}}{(z_{k+1} - z_{k-1})(z_{k+1} - z_k)} \frac{[w(t_{j+1})]^{-1}}{\Gamma(3-\alpha)} C_k^j\,, \tag{50}$$

$$A_k^j = (2-\alpha)\big(2z_{j+1} - z_k - z_{k+1}\big)p_k^j - 2(1-\alpha)q_k^j \ , \tag{51}$$

$$B_k^j \ = (2-\alpha)\big(2z_{j+1} - z_{k-1} - z_{k+1}\big)p_k^j - 2(1-\alpha)q_k^j, \tag{52}$$

$$C_k^j = (2-\alpha)\big(2z_{j+1} - z_{k-1} - z_k\big)p_k^j - 2(1-\alpha)q_k^j \ , \tag{53}$$

**Proof: (i)** Since $A_k^j = \left[(2-\alpha)\big(2z_{j+1} - z_k - z_{k+1}\big)p_k^j - 2(1-\alpha)q_k^j\right]$

From the Eq. (13) and Eq. (14), we have

$$A_k^j = (2-\alpha)\big(2z_{j+1} - z_k - z_{k+1}\big)\left[\big(z_{j+1} - z_k\big)^{1-\alpha} - \big(z_{j+1} - z_{k+1}\big)^{1-\alpha}\right]$$

$$-2(1-\alpha)\left[\big(z_{j+1} - z_k\big)^{2-\alpha} - \big(z_{j+1} - z_{k+1}\big)^{2-\alpha}\right]$$

$$= (2-\alpha)\big(z_{j+1} - z_k + z_{j+1} - z_{k+1}\big)\left[\big(z_{j+1} - z_k\big)^{1-\alpha} - \big(z_{j+1} - z_{k+1}\big)^{1-\alpha}\right]$$

$$-2\left[\big(z_{j+1} - z_k\big)^{2-\alpha} - \big(z_{j+1} - z_{k+1}\big)^{2-\alpha}\right] + 2\alpha\left[\big(z_{j+1} - z_k\big)^{2-\alpha} - \big(z_{j+1} - z_{k+1}\big)^{2-\alpha}\right]$$

$$= (2-\alpha)\left[\big(z_{j+1} - z_k\big)^{2-\alpha} - \big(z_{j+1} - z_{k+1}\big)^{2-\alpha}\right]$$

$$(2-\alpha)\left[\big(z_{j+1} - z_{k+1}\big)\big(z_{j+1} - z_k\big)^{1-\alpha} - \big(z_{j+1} - z_k\big)\big(z_{j+1} - z_{k+1}\big)^{1-\alpha}\right]$$

$$-2\left[\big(z_{j+1} - z_k\big)^{2-\alpha} - \big(z_{j+1} - z_{k+1}\big)^{2-\alpha}\right] + 2\alpha\left[\big(z_{j+1} - z_k\big)^{2-\alpha} - \big(z_{j+1} - z_{k+1}\big)^{2-\alpha}\right]$$

$$= \alpha\left[\big(z_{j+1} - z_k\big)^{2-\alpha} - \big(z_{j+1} - z_{k+1}\big)^{2-\alpha}\right]$$

$$-(2-\alpha)\left[\big(z_{j+1} - z_k\big)\big(z_{j+1} - z_{k+1}\big)^{1-\alpha} - \big(z_{j+1} - z_{k+1}\big)\big(z_{j+1} - z_k\big)^{1-\alpha}\right] \tag{54}$$

**(ii)** Since $2\left[\big(z_{j+1} - z_k\big)^{2-\alpha} - \big(z_{j+1} - z_{k+1}\big)^{2-\alpha}\right]$

$$-(2-\alpha)(z_{k+1} - z_k)\left[\big(z_{j+1} - z_k\big)^{1-\alpha} + \big(z_{j+1} - z_{k+1}\big)^{1-\alpha}\right] \tag{55}$$

$= 2\left[\big(z_{j+1} - z_k\big)^{2-\alpha} - \big(z_{j+1} - z_{k+1}\big)^{2-\alpha}\right] - (2-\alpha)\big(z_{j+1} - z_{j+1} + z_{k+1} - z_k\big)\left[\big(z_{j+1} - z_k\big)^{1-\alpha} + \big(z_{j+1} - z_{k+1}\big)^{1-\alpha}\right]$

$= 2\left[\big(z_{j+1} - z_k\big)^{2-\alpha} - \big(z_{j+1} - z_{k+1}\big)^{2-\alpha}\right] - (2-\alpha)\big(z_{j+1} - z_k\big)\left[\big(z_{j+1} - z_k\big)^{1-\alpha} + \big(z_{j+1} - z_{k+1}\big)^{1-\alpha}\right] + (2-\alpha)\big(z_{j+1} - z_{k+1}\big)\left[\big(z_{j+1} - z_k\big)^{1-\alpha} + \big(z_{j+1} - z_{k+1}\big)^{1-\alpha}\right]$

$= 2\big(z_{j+1} - z_k\big)^{2-\alpha} - 2\big(z_{j+1} - z_{k+1}\big)^{2-\alpha} - (2-\alpha)\big(z_{j+1} - z_k\big)^{2-\alpha} - (2-\alpha)\big(z_{j+1} - z_k\big)\big(z_{j+1} - z_{k+1}\big)^{1-\alpha} + (2-\alpha)\big(z_{j+1} - z_{k+1}\big)\big(z_{j+1} - z_k\big)^{1-\alpha} + (2-\alpha)\big(z_{j+1} - z_{k+1}\big)^{2-\alpha}$

$= \alpha\big(z_{j+1} - z_k\big)^{2-\alpha} - \alpha\big(z_{j+1} - z_{k+1}\big)^{2-\alpha} - (2-\alpha)\big(z_{j+1} - z_k\big)\big(z_{j+1} - z_{k+1}\big)^{1-\alpha} + (2-\alpha)\big(z_{j+1} - z_{k+1}\big)\big(z_{j+1} - z_k\big)^{1-\alpha}$

$= \alpha\left[\big(z_{j+1} - z_k\big)^{2-\alpha} - \big(z_{j+1} - z_{k+1}\big)^{2-\alpha}\right] - (2-\alpha)\big(z_{j+1} - z_k\big)\big(z_{j+1} - z_{k+1}\big)^{1-\alpha} + (2-\alpha)\big(z_{j+1} - z_{k+1}\big)\big(z_{j+1} - z_k\big)^{1-\alpha}$

$= \alpha\left[\big(z_{j+1} - z_k\big)^{2-\alpha} - \big(z_{j+1} - z_{k+1}\big)^{2-\alpha}\right] - (2-\alpha)\left\{\big(z_{j+1} - z_k\big)\big(z_{j+1} - z_{k+1}\big)^{1-\alpha} - \big(z_{j+1} - z_{k+1}\big)\big(z_{j+1} - z_k\big)^{1-\alpha}\right\}$

$= A_k^j$ (56)

**(iii)** Since $\big(z_{j+1} - z_{k+1}\big) < \big(z_{j+1} - z_k\big)$

$\Rightarrow \big(z_{j+1} - z_{k+1}\big)^{1-\alpha} < \big(z_{j+1} - z_k\big)^{1-\alpha}$

$\Rightarrow \big(2z_{j+1} - z_{k+1} - z_k\big)\big(z_{j+1} - z_{k+1}\big)^{1-\alpha} < \big(2z_{j+1} - z_{k+1} - z_k\big)\big(z_{j+1} - z_k\big)^{1-\alpha}$

$\Rightarrow \big(z_{j+1} - z_{k+1} + z_{j+1} - z_k\big)\big(z_{j+1} - z_{k+1}\big)^{1-\alpha} < \big(z_{j+1} - z_{k+1} + z_{j+1} - z_k\big)\big(z_{j+1} - z_k\big)^{1-\alpha}$

$\Rightarrow \big(z_{j+1} - z_{k+1}\big)^{2-\alpha} + \big(z_{j+1} - z_k\big)\big(z_{j+1} - z_{k+1}\big)^{1-\alpha} < \big(z_{j+1} - z_{k+1}\big)\big(z_{j+1} - z_k\big)^{1-\alpha} + \big(z_{j+1} - z_k\big)^{2-\alpha}$

$\Rightarrow \big(z_{j+1} - z_{k+1}\big)^{2-\alpha} - \big(z_{j+1} - z_k\big)^{2-\alpha} < \big(z_{j+1} - z_{k+1}\big)\big(z_{j+1} - z_k\big)^{1-\alpha} - \big(z_{j+1} - z_k\big)\big(z_{j+1} - z_{k+1}\big)^{1-\alpha}$ (57)

$\Rightarrow \alpha\left\{\big(z_{j+1} - z_{k+1}\big)^{2-\alpha} - \big(z_{j+1} - z_k\big)^{2-\alpha}\right\} < (2-\alpha)\left\{\big(z_{j+1} - z_{k+1}\big)\big(z_{j+1} - z_k\big)^{1-\alpha} - \big(z_{j+1} - z_k\big)\big(z_{j+1} - z_{k+1}\big)^{1-\alpha}\right\}$ (58)

$\Rightarrow \alpha\left\{\big(z_{j+1} - z_k\big)^{2-\alpha} - \big(z_{j+1} - z_{k+1}\big)^{2-\alpha}\right\} - (2-\alpha)\left\{\big(z_{j+1} - z_k\big)\big(z_{j+1} - z_{k+1}\big)^{1-\alpha} - \big(z_{j+1} - z_{k+1}\big)\big(z_{j+1} - z_k\big)^{1-\alpha}\right\} > 0$ (59)

$\Rightarrow A_k^j > 0$

**(iv)** Since $B_k^j - A_k^j = (2-\alpha)\big(2z_{j+1} - z_{k-1} - z_{k+1}\big)p_k^j - 2(1-\alpha)q_k^j$

$-(2-\alpha)\big(2z_{j+1} - z_k - z_{k+1}\big)p_k^j - 2(1-\alpha)q_k^j$

$= (2-\alpha)(z_k - z_{k-1})p_k^j > 0$ (60)

Hence, from Eq. (60), $B_k^j > A_k^j$.

Similarly,

$$C_k^j - B_k^j = (2-\alpha)(z_{k+1} - z_k)p_k^j > 0. \tag{61}$$

Hence, from Eq. (61), $C_k^j > B_k^j$ and finally we get $C_k^j > B_k^j > A_k^j > 0$.

**(v)** According to assumption the weight function $w(t)$ is positive and the scale function $z(t)$ is nonnegative and strictly monotone increasing. Hence, from Eqs. (48), (49) and (50) that $a_k^j$, $b_k^j$ and $c_k^j$ are positive.

## 4. Stability and Convergence Analysis

To establish the convergence analysis through stability of the numerical scheme given in Eq. (22), the Lax-Richertmyer theorem [36] is used.

For the stability of the numerical scheme (22), we rewrite the Eq. (23) as:

$$K_{j+1}U_{j+1} = \left(c_{j-1}^j - b_j^j\right).U_j + V_j \text{ ,} 1 \leq j \leq M-1, \tag{62}$$

where

$$U_j = \left[u_j^1, u_j^2, \cdots, u_j^i, \cdots, u_j^{N-1}\right]^T, \tag{63}$$

$$V_j = \left[v_j^1, v_j^2, \cdots, v_j^i, \cdots, v_j^{N-1}\right]^T, \tag{64}$$

and

$$v_j^i = \left(p_j u_1^i - q_j u_0^i\right) + a_{j-1}^j u_{j-2}^i + \left(a_j^j - b_{j-1}^j\right)u_{j-1}^i + \sum_{k=1}^{j-2}\left(a_k^j u_{k-1}^i - b_k^j u_k^i + c_k^j u_{k+1}^i\right) - f_{j+1}^i \tag{65}$$

For the stability of numerical scheme given in Eq. (23), we have the following stability theorem.

**Theorem 1:** If coefficient $c_j^j$ in obtained matrix $K_{j+1}$ satisfies the monotonic increasing property, then the numerical scheme given in Eq. (23) is convergent as well as stable.

**Proof.** Since the coefficient matrix $K_{j+1}$ is tridiagonal and strictly diagonally dominant hence it is invertible. Therefore Eq. (23) is solvable and iteration scheme is well-posed.

Now, we wish to show the boundedness of posteriori error. Suppose $u_j$ represents the exact solution at $t = t_j$ and $\epsilon_j$ denotes the posteriori error that is $\epsilon_j = u_j - U_j$. Then, we have

$$\epsilon_{j+1} = K_{J+1}^{-1}\left(c_{j-1}^j - b_j^j\right)\epsilon_j + O((\Delta t)^{3-\alpha} + \Delta x^2), \text{ for } 1 \leq j \leq M-1. \tag{66}$$

By [37, Theorem A], suppose $K_{j+1} = \left[k_{ni,nj}\right]_{(N-1)\times(N-1)}$, then

$$\left\|K_{j+1}^{-1}\right\| \leq \frac{1}{\min\limits_{1\leq ni\leq N-1}\left\{\left|k_{ni,nj}\right|-\sum_{nj=1,\ nj\neq ni}^{N-1}\left|k_{ni,nj}\right|\right\}} = \frac{1}{c_j^j}.$$

According to assumption given in Eq. (66), we have

$$\left\|K_{j+1}^{-1}\left(c_{j-1}^{j} - b_j^j\right)\right\| \leq 1. \tag{67}$$

Hence, Eq. (66) and Eq. (67) suggest the existence of a positive constant $C$, such that

$$\|\epsilon_{j+1}\| \leq \left\|K_{j+1}^{-1}\right\|\left|c_{j-1}^{j} - b_j^j\right|.\|\epsilon_j\| + C(\Delta t^{3-\alpha} + \Delta x^2)$$

$$\leq \|\epsilon_j\| + C(\Delta t^{3-\alpha} + \Delta x^2),$$

which implies

$$\|\epsilon_{j+1}\| \leq \left(\|E_j\| + \|E_j\|^2 + \cdots + \|E_j\|^j\right) C(\Delta t^{3-\alpha} + \Delta x^2),$$

where $E_j = K_{j+1}^{-1}\left(c_{j-1}^{j} - b_j^j\right)$ is an amplification matrix. By using stability condition $\|E_j\| \leq 1$ given in Eq. (67) we have

$$\|\epsilon_{j+1}\| \leq jC(\Delta t^{3-\alpha} + \Delta x^2), \quad 1 \leq j \leq M-1,$$

which shows that the error is bounded, and $\|\epsilon_{j+1}\| \to 0$ as $\Delta t \to 0, \Delta x \to 0$. Hence, the numerical scheme given in Eq. (23) is stable and consistent as well. Hence, by the Lax-Richertmyer theorem [36], the convergence of numerical scheme is established.

## 5. Numerical Results and Their Analysis

In the present Section, the validity and numerical ability of our high order proposed finite difference numerical scheme for GFDEs are demonstrated via five examples. And also, the numerical results are compared with the known results in literature. In the first and fourth example, we solve with non-zero initial condition. Although, the other examples are solved by zero initial condition. Fifth example is solved by non-zero boundary condition while others are solved by zero boundary condition. The GFD is defined using the scale and weight function so we also present solution of example 4 with non-zero weight function. Here, we present the numerical results, exact results through figures and tables. The maximum absolute error (MAE) and convergence order (CO) are also presented through tables. Spatial as well as temporal direction errors of the numerical example are also presented through tables.

To check the CO of approximation of GFD, we take $u(x,t) = t - t^3$. It is also clear from Table 1, 2 that the CO in the spatial direction is $(\Delta t)^{3-\alpha}$. The approximation results are obtained for weight function $w(t) = 1$ and it validates the theoretical results presented earlier. To calculate the MAE, and CO , we use the formula, CO =lg[MAE($\Delta t$)/MAE($\Delta t/2$)]/lg(2) , where MAE is obtained using maximum of absolute value of difference of exact solution and numerical solution at the node point.

Table 1: MAE and CO in approximation of GFD for $u(x,t) = t - t^3$ with $z(t) = t^2, \alpha = 0.5$, at $t = 0.6$ .

| $\Delta t$ | MAE | CO |
|---|---|---|
| 1/10 | $6.84003 \times 10^{-3}$ | |
| 1/20 | $1.04274 \times 10^{-3}$ | 2.71362 |
| 1/40 | $1.70108 \times 10^{-4}$ | 2.61586 |
| 1/80 | $2.85787 \times 10^{-5}$ | 2.57344 |
| 1/160 | $4.87585 \times 10^{-6}$ | 2.55121 |

Table 2: MAE and CO in approximation of GFD for $u(x,t) = t - t^3$ with $z(t) = t^2, \alpha = 0.2$, at $t = 0.6$.

| $\Delta t$ | MAE | CO |
|---|---|---|
| 1/10 | $1.48829 \times 10^{-3}$ | |
| 1/20 | $2.07611 \times 10^{-4}$ | 2.8417 |
| 1/40 | $2.99599 \times 10^{-5}$ | 2.79278 |
| 1/80 | $4.35344 \times 10^{-6}$ | 2.78281 |
| 1/160 | $6.31831 \times 10^{-7}$ | 2.78455 |

**Example 1:** Consider GFDEs [30] defined by Eq. (1) as,

$$^*\partial_t^\alpha u(x\,,t) - \frac{\partial^2 u(x\,,t)}{\partial^2 x} = \frac{2x(x-1)t^{2-\alpha}}{\Gamma(3-\alpha)} - 2t^2 + \pi^2 \sin(\pi x),\;\; x \in [0\,,1],\;\; t > 0, \tag{68}$$

with the initial and boundary conditions as, $u(x\,,0) = \sin(\pi x)$, $u(0\,,t) = u(1\,,t) = 0$. We take $z(t) = t,\; w(t) = 1$ and $0 \le t \le 1$.

The exact solution of the Eq. (1) under above condition is given by $u(x\,,t) = x(x-1)t^2 + \sin(\pi x)$. We solved Example 1 with finite difference scheme proposed in Section 2 for different $\alpha = 0.6, 0.4$ and $0.7$, and different step sizes $\Delta x = 1/500, \Delta t = 1/200$; $\Delta x = 1/400, \Delta t = 0.01$, and $\Delta x = 0.05, \Delta t = 0.05$. The numerical results are shown in the Figures 1(a), 1(b) 1(c) and 1(d). Figure 1(a) displays for exact solution while Figures 1(b) 1(c) and 1(d) present for numerical solutions. We also solved Example 1 with 1) varying both $\Delta x$ and $\Delta t$ with $\alpha = 0.85$; 2) varying both $\Delta x$ and $\Delta t$ with $\alpha = 0.6$; 3) varying $\Delta x$ and fixed $\Delta t = 1/600$ with $\alpha = 0.85$ and 4) varying $\Delta t$ and fixed $\Delta x = 1/512$ with $\alpha = 0.85$ which are shown through Table1 (a), 1(b), 1(c) and 1(d) respectively. It is clear from the Tables 1(a) and 1(b) that the finite difference scheme is of second order of convergence. Table 1(c) also shows the numerical scheme having second order of convergence in spatial direction. Figure 1 shows that whenever we reduce the step size in spatial or temporal or both directions, the numerical results converge to the exact solution which validates the scheme is numerically convergent. Table 1(d) shows that MAE and CO for $\alpha = 0.85$ and compares the results with reference [30] for Example 1. It is clear from the Table 1(c) and Table 1(d) our proposed numerical scheme is of high order of convergence and achieve the good accuracy. The numerical scheme proposed in [30] is $O(\Delta x^2 + \Delta t)$ while our numerical finite difference scheme of order $O(\Delta x^2 + \Delta t^{3-\alpha})$.

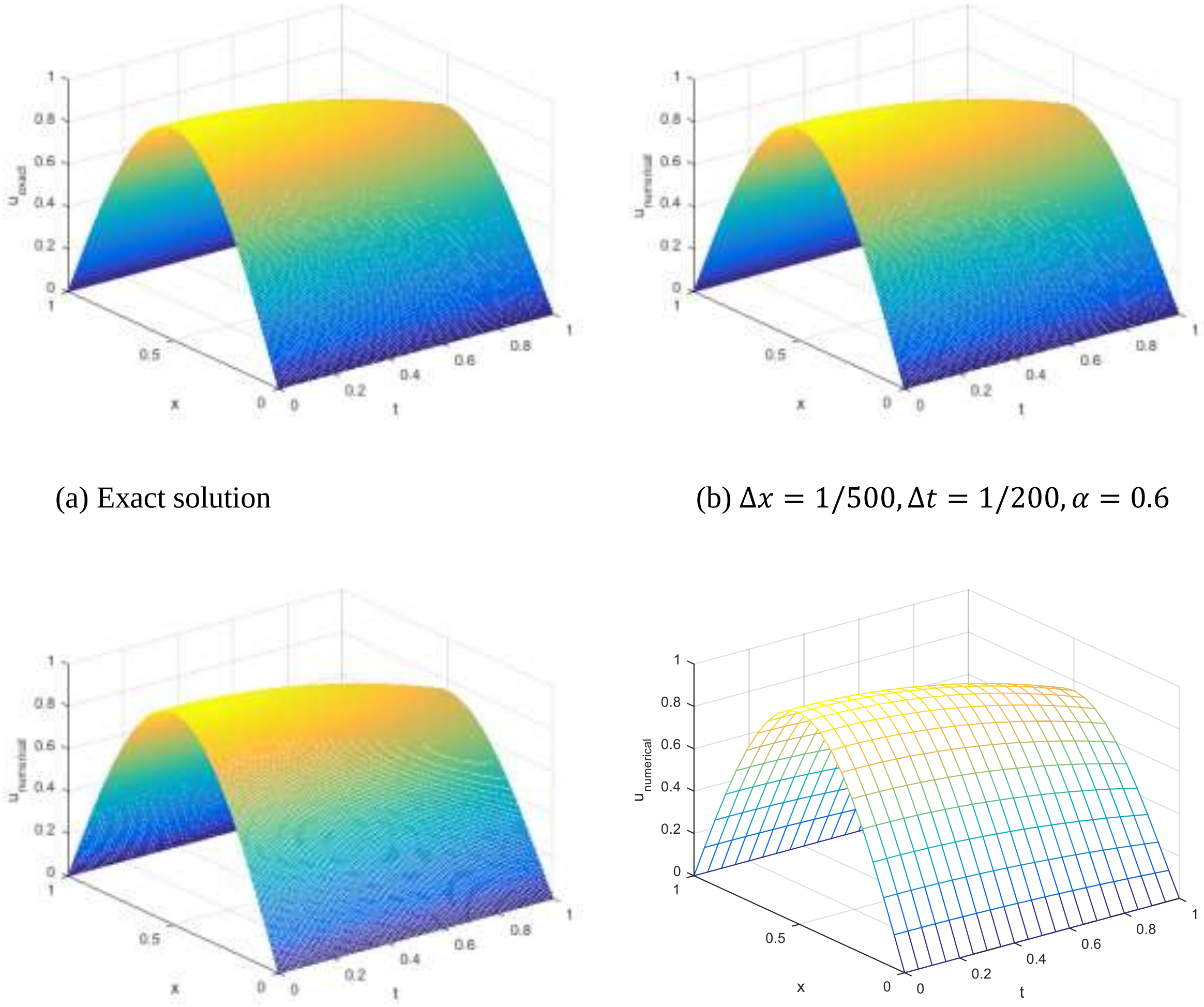

(a) Exact solution

(b) $\Delta x = 1/500, \Delta t = 1/200, \alpha = 0.6$

(c) $\Delta x = 1/400, \Delta t = 0.01, \alpha = 0.4$

(d) $\Delta x = 0.05, \Delta t = 0.05, \alpha = 0.7$

**Fig 1:** Comparison of solutions for Example 1 with different parameters.

Table 1(a): MAE and CO for Example 1 with $\alpha = 0.85,\ \delta = 1\ z(t) = t,\ w(t) = 1$.

| $\Delta t$ | $\Delta x$ | MAE | CO |
|---|---|---|---|
| 1/8 | 1/8 | 0.012682479250020 | - |
| 1/16 | 1/16 | 0.003154209795232 | 2.0075 |
| 1/32 | 1/32 | 0.000787641393682 | 2.0017 |
| 1/64 | 1/64 | 0.000196866793636 | 2.0003 |
| 1/128 | 1/128 | 0.000049216024556 | 2.0000 |

Table 1(b): MAE and CO for Example 1 with $\alpha = 0.6,\ \delta = 1\ z(t) = t,\ w(t) = 1$.

| $\Delta t$ | $\Delta x$ | MAE | CO |
|---|---|---|---|
| 1/8 | 1/8 | 0.012307702788089 | - |
| 1/16 | 1/16 | 0.003063419371486 | 2.0063 |
| 1/32 | 1/32 | 0.000765192639171 | 2.0012 |
| 1/64 | 1/64 | 0.000191279519140 | 2.0001 |
| 1/128 | 1/128 | 0.000047821744845 | 1.9999 |

Table 1(c): MAE and CO for Example 1 with $\alpha = 0.85,\ \delta = 1\ z(t) = t,\ w(t) = 1$, and $\Delta t = 1/600$.

| $\Delta x$ | MAE [30] | CO [30] | MAE | CO |
|---|---|---|---|---|
| 1/8 | 0.05290053759936 | - | 0.012697453071990 | - |
| 1/16 | 0.01365644886944 | 1.9537 | 0.003156727061421 | 2.0080 |
| 1/32 | 0.00384054346854 | 1.8302 | 0.000788084897559 | 2.0020 |
| 1/64 | 0.00108689342074 | 1.8211 | 0.000196952333293 | 2.0005 |

Table 1(d): MAE and CO for Example 1 for $\alpha = 0.85$ with different $\Delta t$ and fixing $\Delta x = 1/512$.

| $\Delta t$ | MAE [30] | CO [30] | MAE | CO |
|---|---|---|---|---|
| 1/8 | 0.06935363245950 | - | 0.001149644902625 | - |
| 1/16 | 0.03256819563577 | 1.0905 | 0.000393504352298 | 1.5467 |
| 1/32 | 0.01578117824480 | 1.0453 | 0.000123282907375 | 1.6744 |
| 1/64 | 0.00777362866697 | 1.0215 | 0.000035575711316 | 1.7930 |

**Example 2:** Consider GFDEs defined in Eq. (1) as,

$$ {}^{*}\partial_t^{\alpha} u(x\,,t) - \frac{\partial^2 u(x\,,t)}{\partial^2 x} = \frac{2x(x-1)t^{2-\alpha}}{\Gamma(3-\alpha)} - 2t^2,\; x \in [0\,,1],\; t > 0, \tag{31} $$

with the initial and boundary conditions as $u(x\,,0) = 0$, $u(0\,,t) = u(1\,,t) = 0$. We take $z(t) = t$, $w(t) = 1$ and $0 \le t \le 1$.

The exact solution of the Eq. (1) under above conditions is given by $u(x\,,t) = x(x-1)t^2$. We solved Example 2 with finite difference scheme proposed in Section 2 for different $\alpha = 0.4, 0.5$ and 0.8 and different step sizes $\Delta x = 1/200, \Delta t = 0.01$; $\Delta x = 1/200, \Delta t = 0.05$ and $\Delta x = 0.05, \Delta t = 0.05$. The numerical results are shown in the Figures 2(a), 2(b) 2(c) and 2(d). The Figure 2(a) display for the exact solution while Figures 2(b) 2(c) and 2(d) present for numerical solution. We also solved Example 2 with 1) varying both $\Delta x$ and $\Delta t$ with $\alpha = 0.8$; and 2) varying $\Delta t$ and fixed $\Delta x = 1/5000$ with $\alpha = 0.5$ which are shown through Table 2 (a) and 2(b) respectively. It is clear from the Table 2(a)) that our finite difference of second order of convergence. Table 2(b) also shows that the numerical scheme having $3 - \alpha$ order of convergence in temporal direction. It is clear from the Table 2(a) and Table 1(b) our proposed numerical scheme is high order of convergence and achieves the good accuracy.

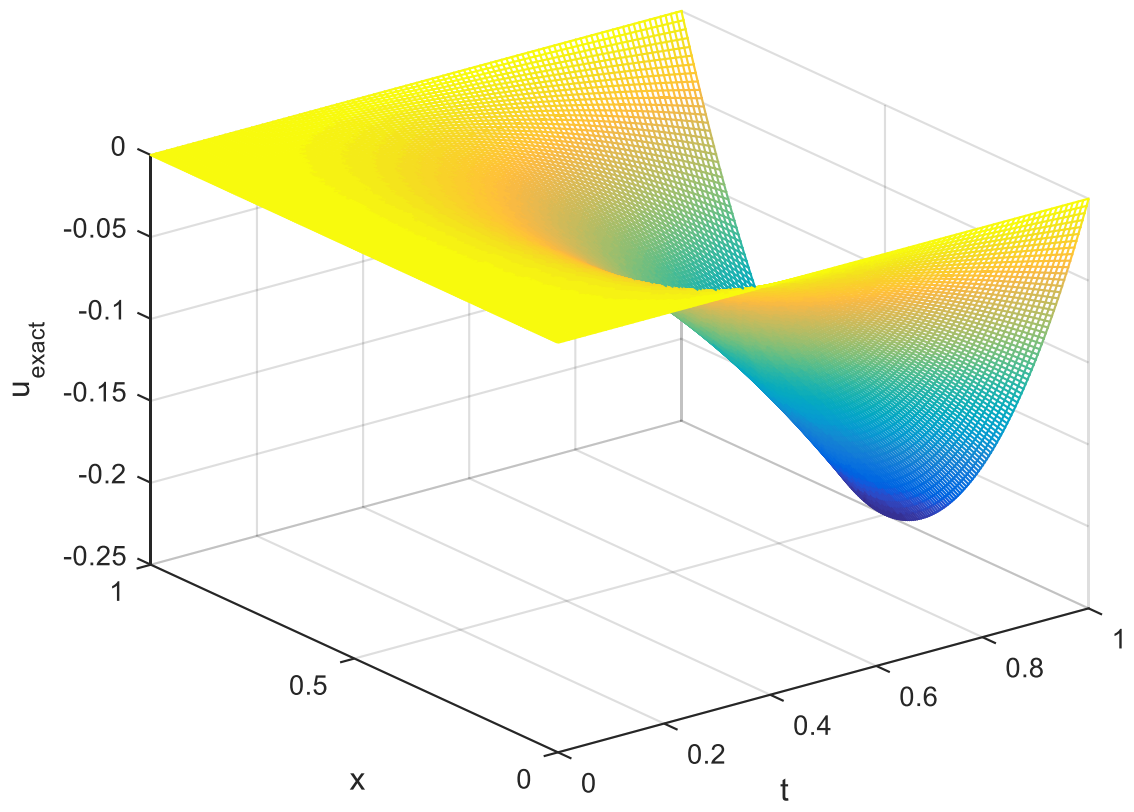

(a) Exact Solution

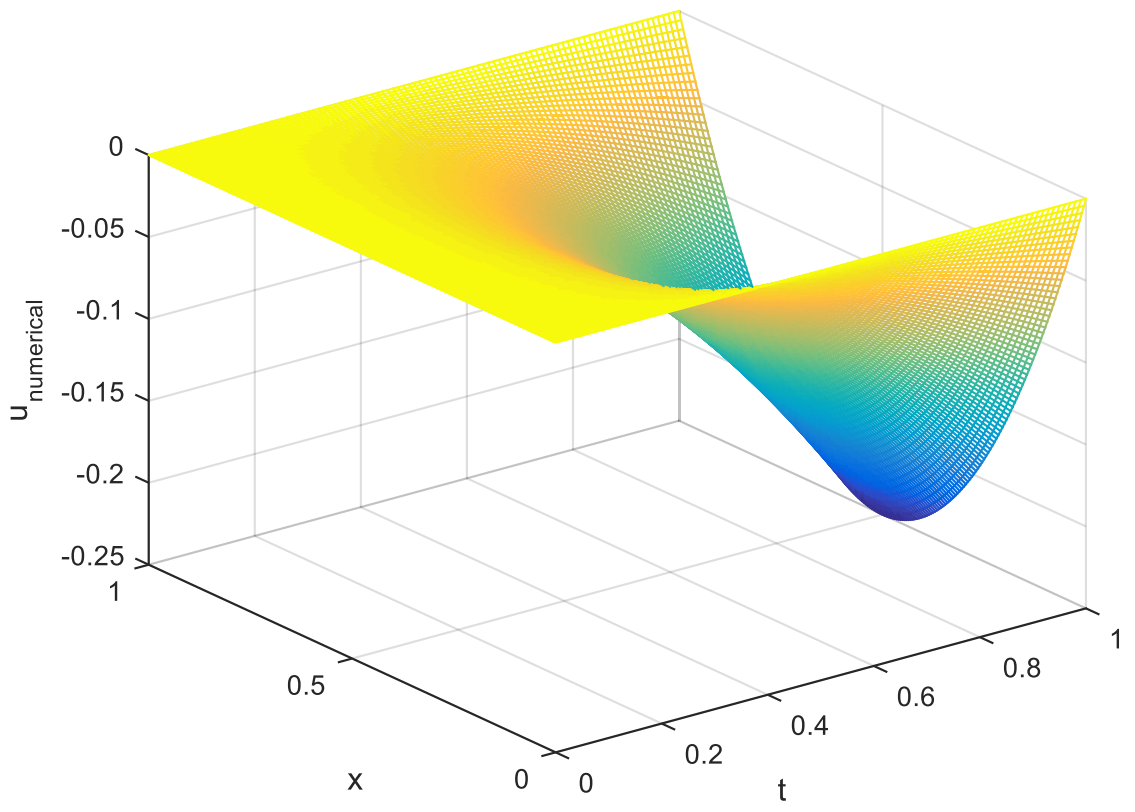

(b) $\Delta x = 1/200, \Delta t = 0.01, \alpha = 0.4$

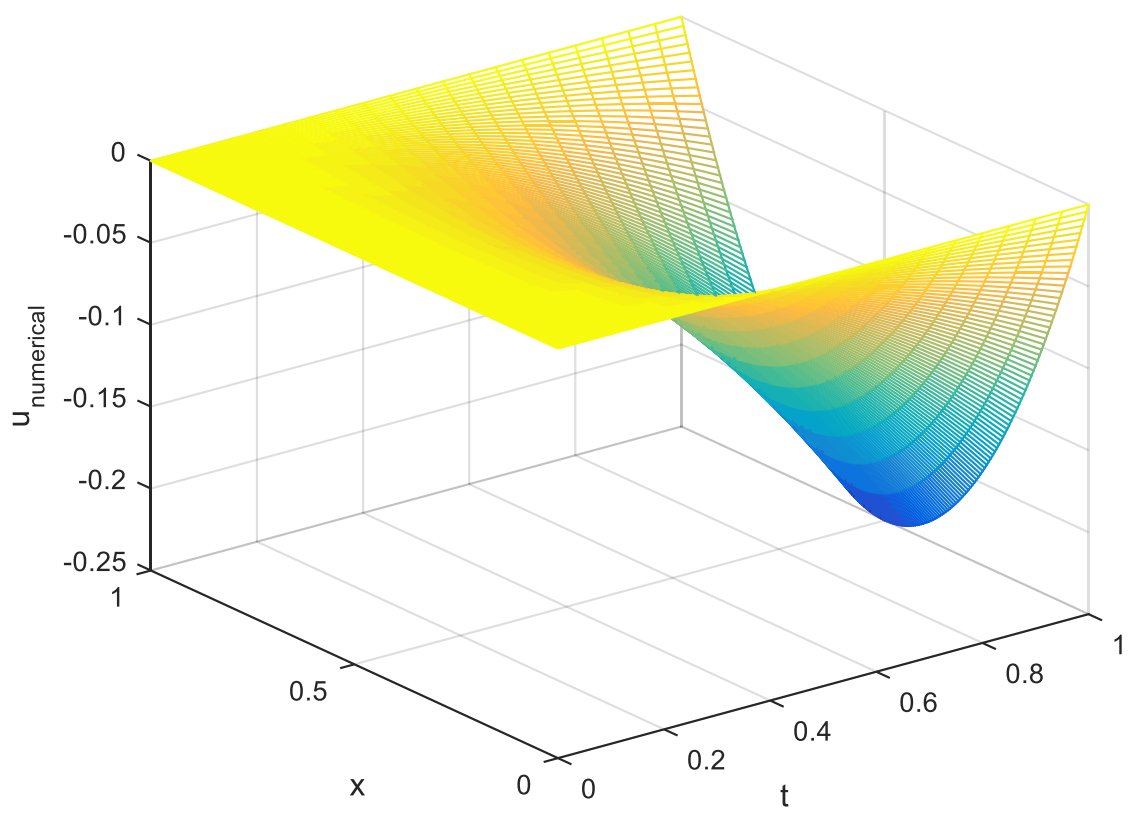

(c) $\Delta x = 1/200, \Delta t = 0.05, \alpha = 0.5$

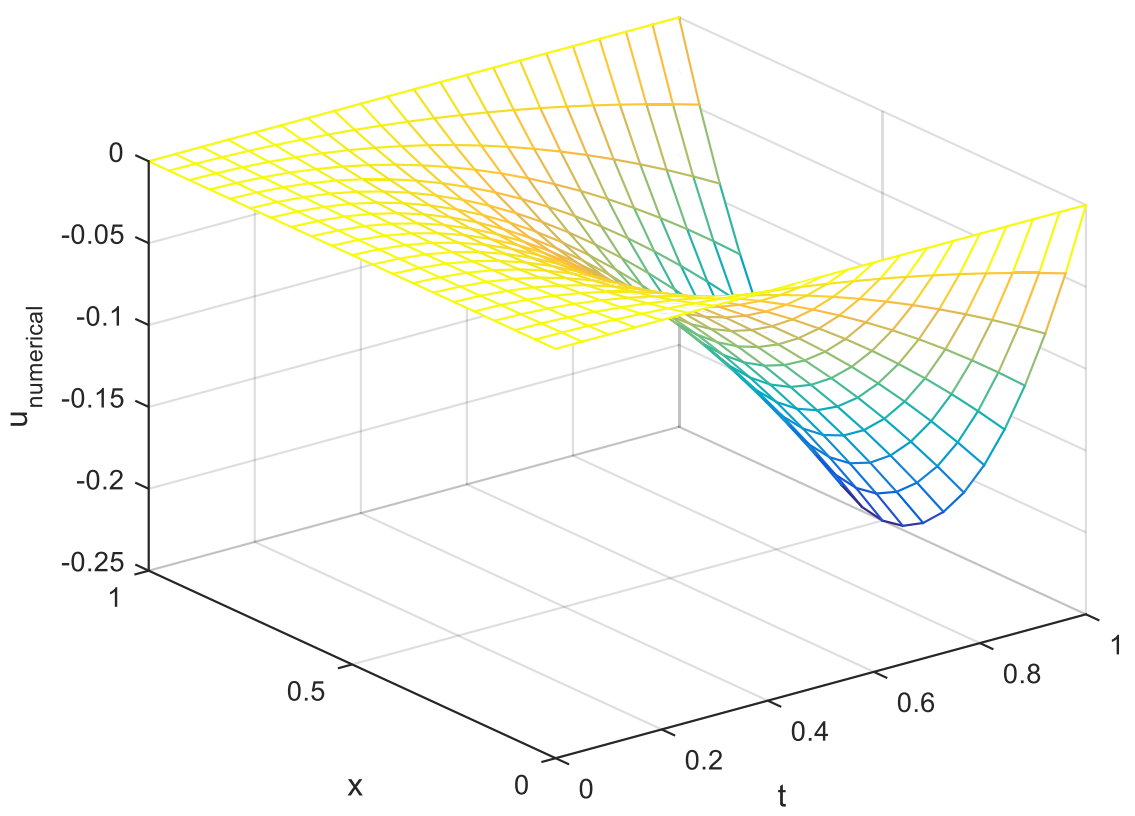

(d) $\Delta x = 0.05, \Delta t = 0.05, \alpha = 0.8$

**Fig 2:** Comparison of solutions of Example 2 with different parameters.

Table 2(a): MAE and CO for Example 2 with $\alpha = 0.8,\ \delta = 1\ z(t) = t,\ w(t) = 1.$

| $\Delta t$ | $\Delta x$ | MAE | CO |
|---|---|---|---|
| 1/10 | 1/10 | 7.0493e-04 | - |
| 1/20 | 1/20 | 2.3408e-04 | 1.5905 |
| 1/40 | 1/40 | 7.2255e-05 | 1.6958 |
| 1/80 | 1/80 | 2.0862e-05 | 1.7922 |
| 1/160 | 1/160 | 5.9576e-06 | 1.8081 |
| 1/320 | 1/320 | 1.6449e-06 | 1.8567 |
| 1/640 | 1/640 | 4.3554e-07 | 1.9172 |

Table 2(b): MAE and CO for Example 2 with $\alpha = 0.5, \Delta x = \frac{1}{5000}, \delta = 1$ for Example 2 for $T = 1$.

| $\Delta t$ | MAE | CO |
|---|---|---|
| 1/10 | 2.6296e-06 | - |
| 1/20 | 3.8045e-07 | 2.7891 |
| 1/40 | 5.8002e-08 | 2.7135 |
| 1/80 | 9.1278e-09 | 2.6678 |
| 1/160 | 1.4762e-09 | 2.6284 |
| 1/320 | 2.7062e-10 | 2.4475 |

**Example 3:** Consider GFDEs defined in Eq. (1) as,

$$ {}^{*}\partial_t^{\alpha} u(x\,,t) - \frac{\partial^2 u(x\,,t)}{\partial^2 x} = \frac{2t^{2-\alpha}}{\Gamma(3-\alpha)} \sin(\pi x) + \pi^2\, t^2 \sin(\pi x)\,,\ x \in [0\,,1],\ t > 0, \qquad (34)$$

with the initial and boundary conditions as $u(x\,,0) = 0, u(0\,,t) = u(1\,,t) = 0$.

The exact solution of the Eq. (1) under above conditions is given by $u(x\,,t) = t^2 \sin(\pi x)$. We solved Example 3 with finite difference scheme proposed in Section 2 with $z(t) = t,\ w(t) = 1$ for different $\alpha = 0.6$ and $0.8$ and different step sizes $\Delta x = 0.001, \Delta t = 0.01$; $\Delta x = 0.01, \Delta t = 0.02$ and $\Delta x = 0.01, \Delta t = 0.025$. The numerical results are shown in the Figures 3(a), 3(b) 3(c) and 3(d). Figure 3(a) displays the exact solution while 3(b) 3(c) and 3(d) present for numerical solutions. We also solved Example 3 with 1) varying both $\Delta x$ and $\Delta t$ with $\alpha = 0.4$; and 2) varying $\Delta x$ and fixed $\Delta t = 1/600$ with $\alpha = 0.4$ which are shown through Table 3 (a) and 3(b) respectively. It is clear from the Tables 3(a) and 3(b) that our finite difference is of second order of convergence. Table 3(b) shows that the numerical scheme having second order of convergence in spatial direction. Figure 3 shows that whenever we reduce the step size in spatial or temporal or both directions, the numerical results convergent to exact solution. which validates the scheme is numerically convergence. Table 3(a) and Table 3(b) are presented for the MAE and CO with $\alpha = 0.4$ . It is clear from the Table 3(a) and Table 3(b) our proposed numerical scheme is high order of convergence and achieve the good accuracy. It is clear from this example that numerical results conclude the theoretical results.

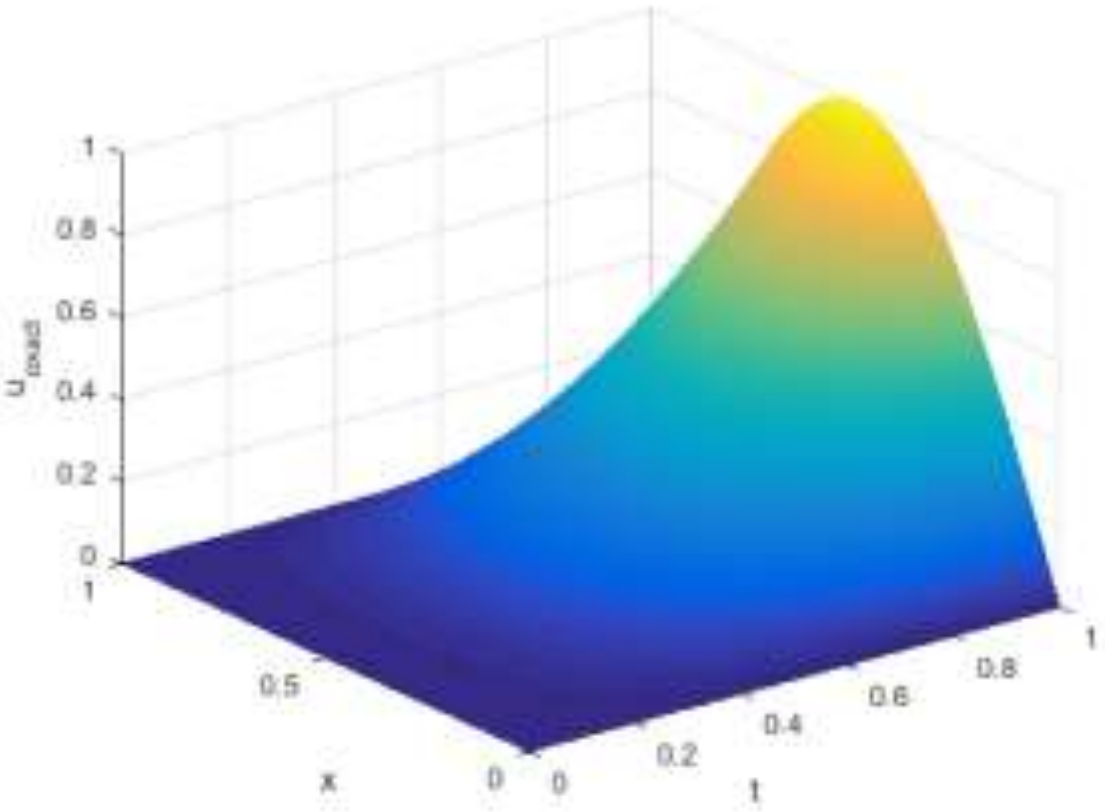

(a) Exact Solution.

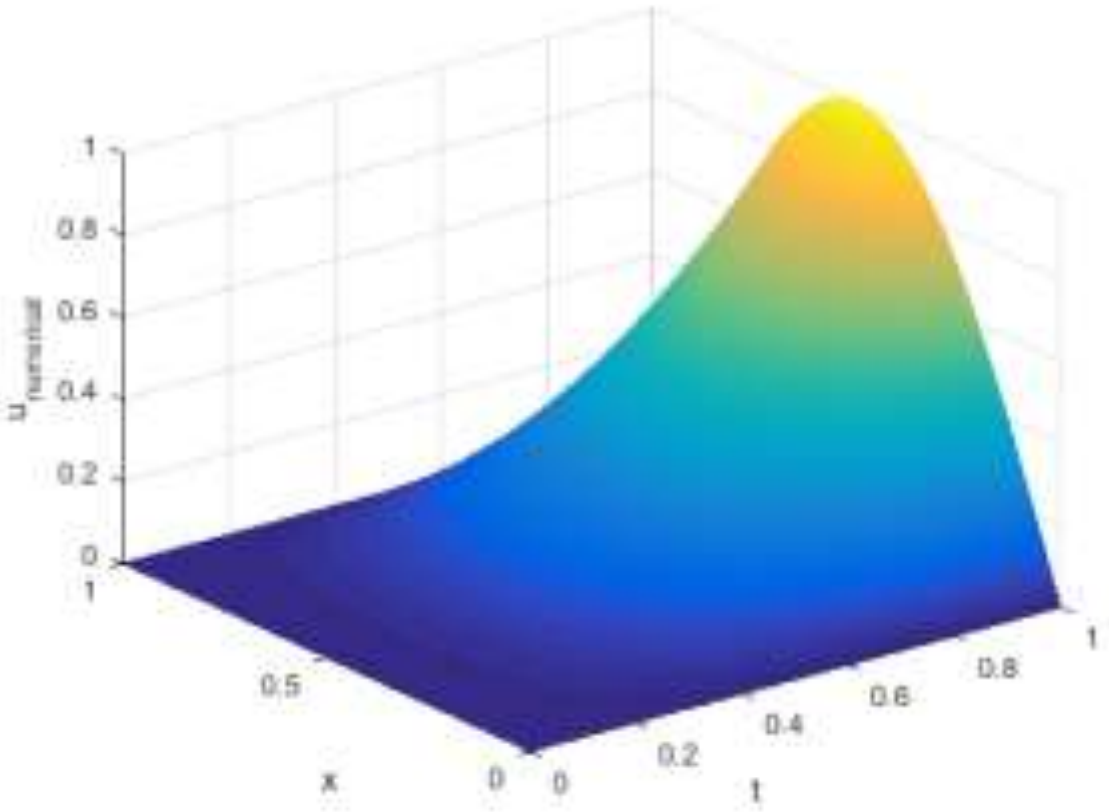

(b) $\Delta x = 0.001, \Delta t = 0.01, \alpha = 0.6.$

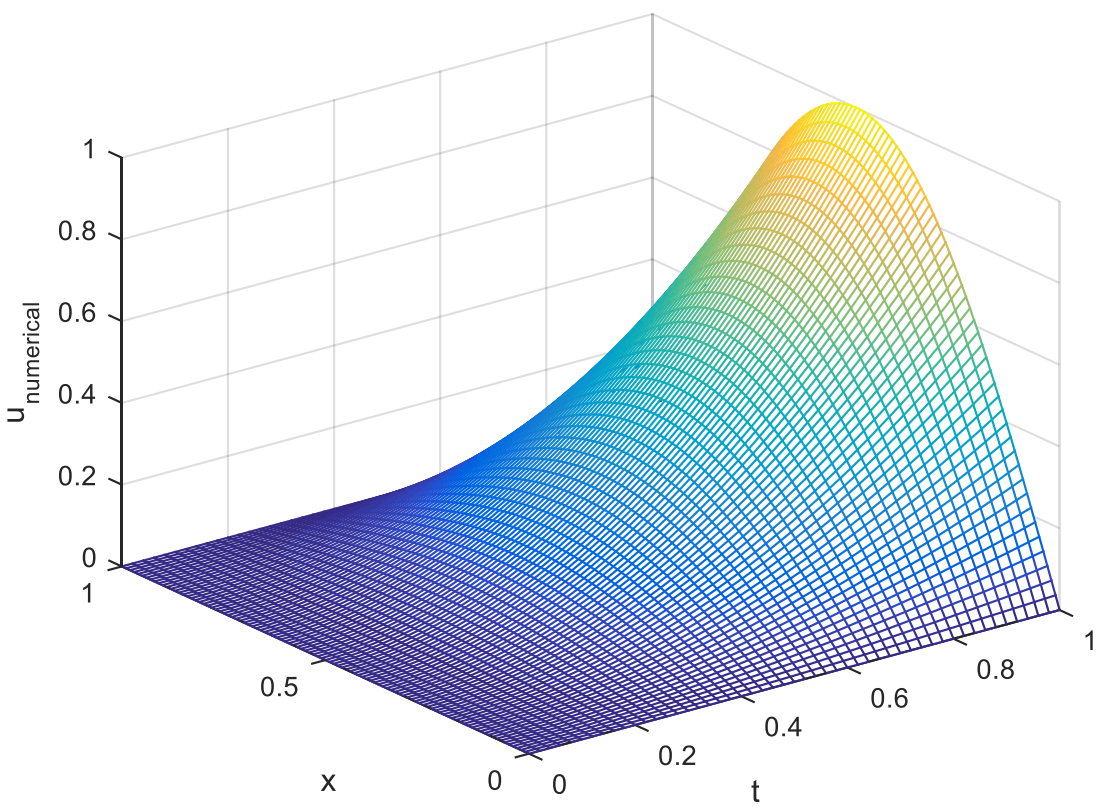

(c) $\Delta x = 0.01, \Delta t = 0.02, \alpha = 0.6.$

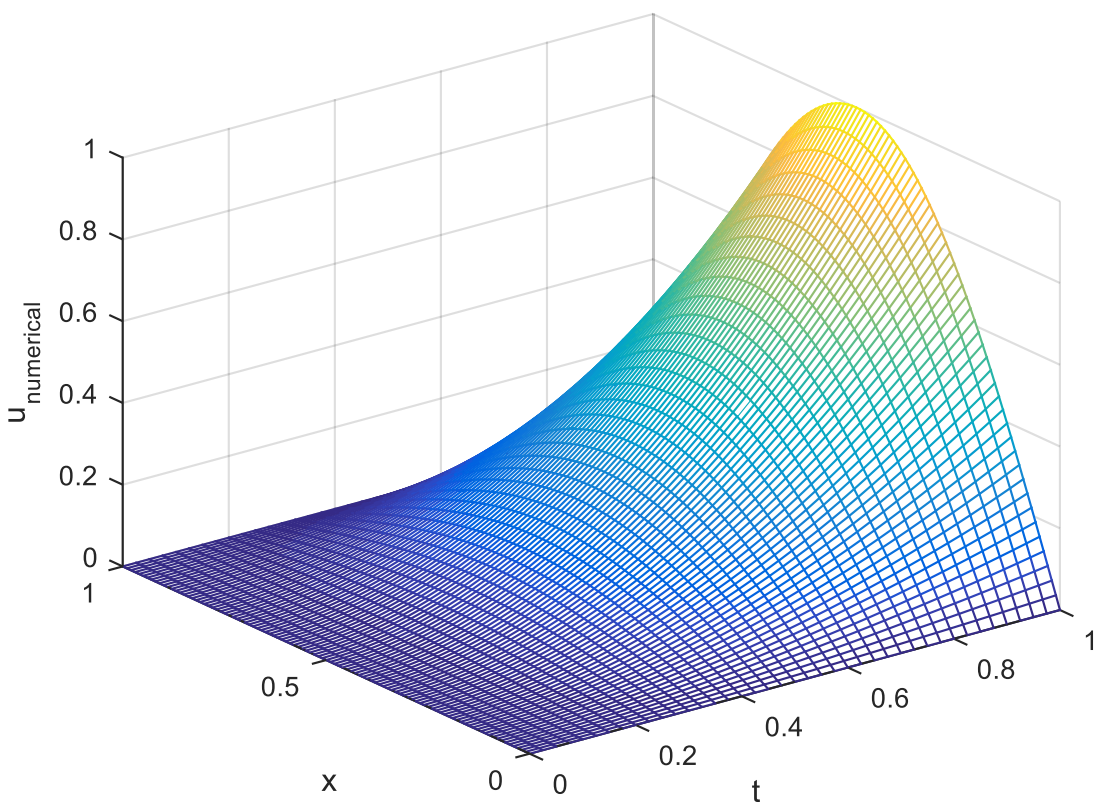

(d) $\Delta x = 0.01, \Delta t = 0.025,\ \alpha = 0.8.$

**Fig 3:** Comparison of solutions of Example 3 with different parameters.

Table 3(a): MAE and CO for Example 3 with $\alpha = 0.4,\ \delta = 1\ z(t) = t,\ w(t) = 1.$

| $\boldsymbol{\Delta t}$ | $\boldsymbol{\Delta x}$ | MAE | CO |
|---|---|---|---|
| 1/10 | 1/10 | 7.2313e-03 | - |
| 1/20 | 1/20 | 1.8017e-03 | 2.0049 |
| 1/40 | 1/40 | 4.4997e-04 | 2.0014 |
| 1/80 | 1/80 | 1.1246e-04 | 2.0005 |
| 1/160 | 1/160 | 2.8111e-05 | 2.0002 |
| 1/320 | 1/320 | 7.0274e-06 | 2.0001 |

Table 3(b): MAE and CO for Example 3 with $\alpha = 0.4,\ \delta = 1\ z(t) = t,\ w(t) = 1\ \Delta t = 1/600$.

| $\Delta x$ | MAE | **CO** |
|---|---|---|
| 1/10 | 7.2237e-03 | - |
| 1/20 | 1.8006e-03 | 2.0042 |
| 1/40 | 4.4983e-04 | 2.0011 |
| 1/80 | 1.1244e-04 | 2.0003 |
| 1/160 | 2.8108e-05 | 2.0001 |
| 1/320 | 7.0270e-06 | 2.0000 |

**Example 4:** Consider GFDEs defined in Eq. (1) as,

$$^{*}\partial_t^{\alpha} u(x\,,t) - \frac{\partial^2 u(x\,,t)}{\partial^2 x} = \frac{1}{\Gamma(1-\alpha)} \sin(\pi x) \exp(-t)\big(\Gamma(1-\alpha) - \gamma(1-\alpha, t)\big) + \pi^2 \sin(\pi x) \exp(-t)\,,$$
$$x \in [0\,,1],\ t > 0, \qquad (35)$$

with the initial and boundary conditions as $u(x\,,0) = \sin(\pi x)$, $u(0\,,t) = u(1\,,t) = 0$. Where $\gamma(s,t)$ denotes the incomplete Gamma function. We take $z(t) = t,\ w(t) = \exp(2t)$ and $0 \le t \le 1$.

The exact solution of the Eq. (1) under above condition is given by $u(x\,,t) = \sin(\pi x)\exp(-t)$. We solved Example 4 with finite difference scheme proposed in Section 2 for $\alpha =\ 0.4$ and different step sizes $\Delta x = \Delta t = 0.001$ and $\Delta x = 0.01, \Delta t = 0.02$. The numerical results are shown in the Figures 4(a), 4(b) and 4(c). The exact solution is displayed through Figure 4(a) while Figures 4(b) and 4(c) show the numerical solutions. We also solved Example 4 with 1) varying both $\Delta x$ and $\Delta t$ with $\alpha = 0.6$; 2) varying $\Delta x$ and fixed $\Delta t = 1/600$ with $\alpha = 0.8$ and 3) varying $\Delta t$ and fixed $\Delta x = 1/10000$ with $\alpha = 0.6$ which are shown through Table 4(a), 4(b) and 4(c) respectively. It is clear from the Tables 4(a), 4(b) and 4(c) that our finite difference of second order of convergence. Table 4(b) and Table 4(c) is also shows that the numerical scheme having second order of convergence in spatial direction and $3 - \alpha$ order of convergence in temporal direction respectively. Table 4(a), 4(b) and 4(c) are presented for the MAE and CO with $\alpha =\ 0.6$ and 0.8. It is clear from the Table 4(a), 4(b) and 4(c) our proposed numerical scheme is high order of convergence and achieves good accuracy in approximating numerical solutions.

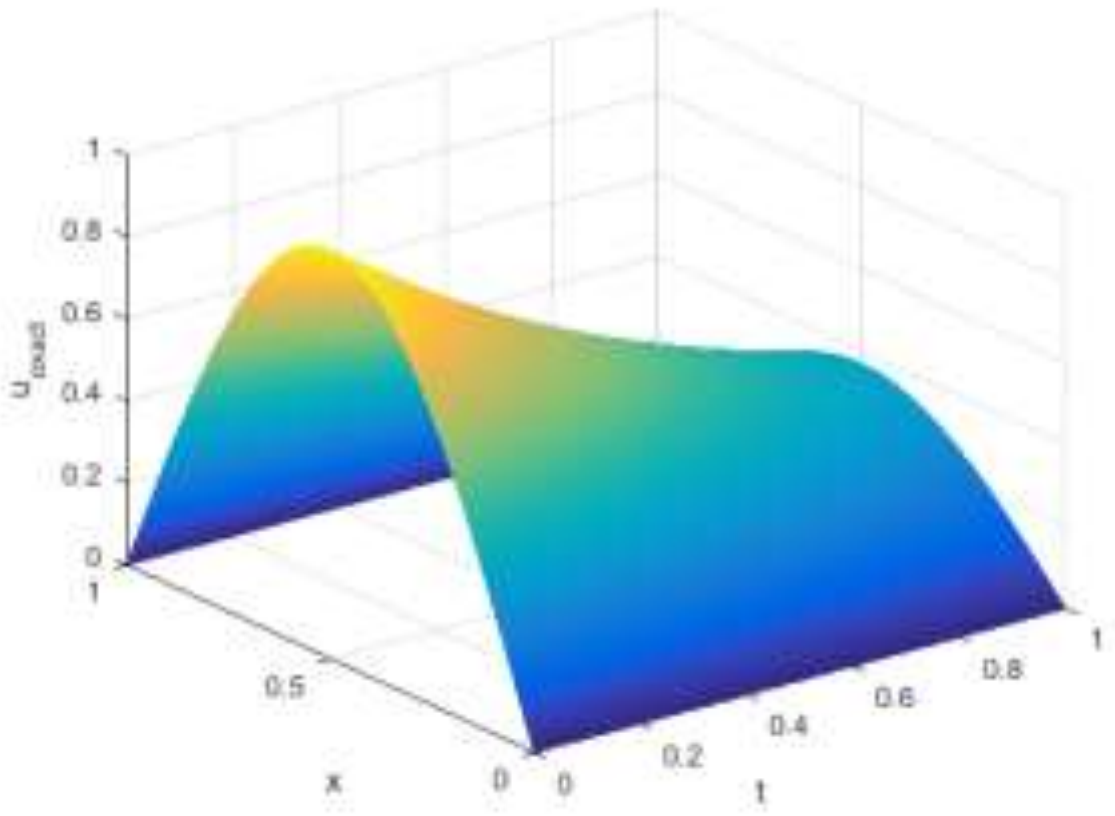

(a) Exact Solution.

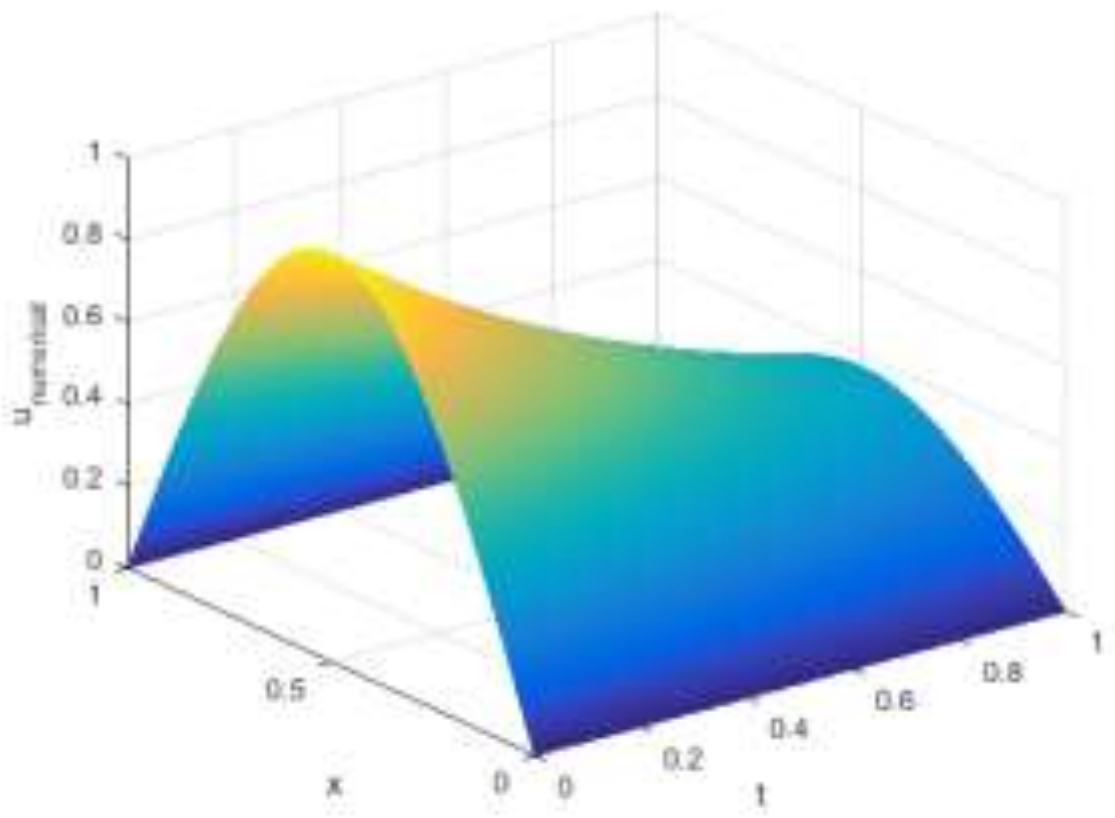

(b) $\Delta x = \Delta t = 0.001, \alpha = 0.4.$

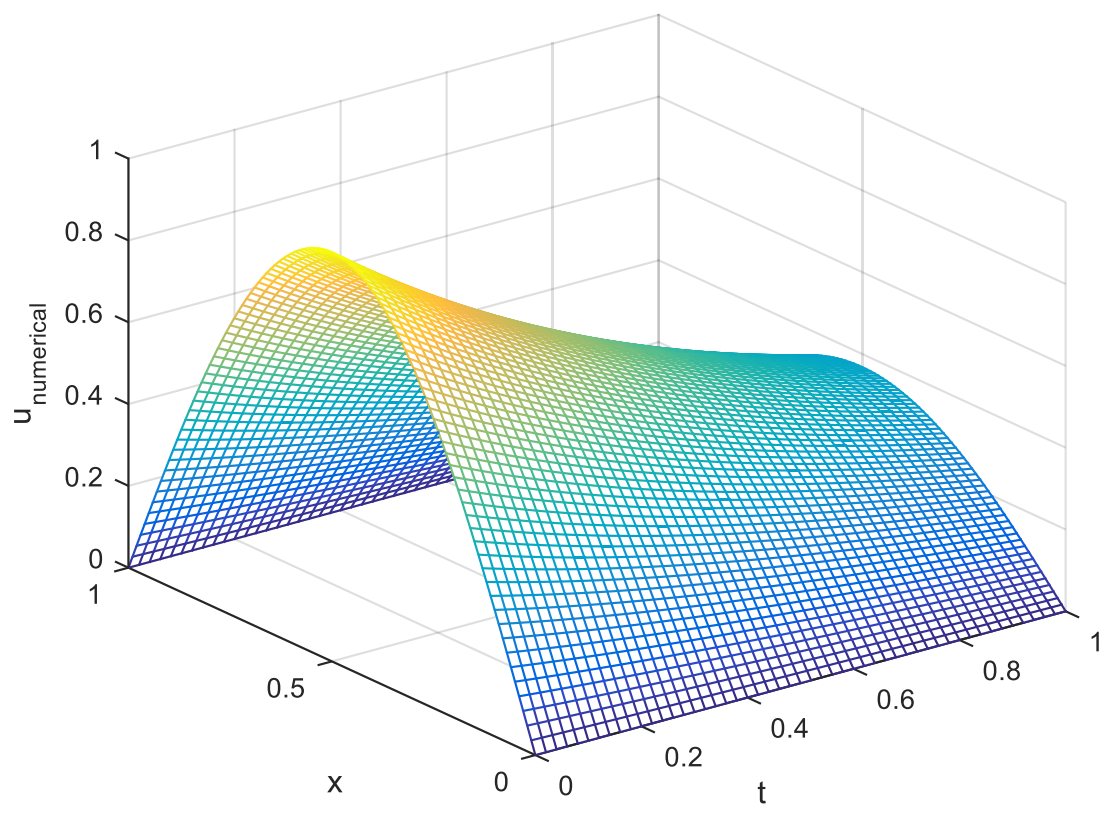

(c) $\Delta x = 0.01, \Delta t = 0.02, \alpha = 0.4.$

**Fig 4:** Comparison of solutions of Example 4 with different parameters.

Table 4(a): MAE and CO for Example 4 with $\alpha = 0.6,\ \delta = 1\ z(t) = t,\ w(t) = \exp(2t).$

| $\Delta t$ | $\Delta x$ | MAE | CO |
|---|---|---|---|
| 1/10 | 1/10 | 5.9274e-03 | - |
| 1/20 | 1/20 | 1.5408e-03 | 1.9437 |
| 1/40 | 1/40 | 3.7698e-04 | 2.0311 |
| 1/80 | 1/80 | 9.3308e-05 | 2.0144 |
| 1/160 | 1/160 | 2.3282e-05 | 2.0028 |
| 1/320 | 1/320 | 5.8117e-06 | 2.0022 |

Table 4(b): MAE and CO for Example 4 with $\alpha = 0.8,\ \delta = 1\ z(t) = t,\ w(t) = \exp(2t)\ \Delta t = 1/600.$

| $\Delta x$ | MAE | CO |
|---|---|---|
| 1/10 | 5.6650e-03 | - |
| 1/20 | 1.4133e-03 | 2.0030 |
| 1/40 | 3.5316e-04 | 2.0006 |
| 1/80 | 8.8310e-05 | 1.9997 |
| 1/160 | 2.2107e-05 | 1.9980 |
| 1/320 | 5.5580e-06 | 1.9919 |

Table 4(c): MAE and CO for Example 4 with $\alpha = 0.6, \Delta x = \frac{1}{10000}, \delta = 1\ z(t) = t,\ w(t) = \exp(2t).$ for $T = 1$

| $\Delta t$ | MAE | CO |
|---|---|---|
| 1/10 | 2.0563e-05 | - |
| 1/20 | 4.0706e-06 | 2.3367 |
| 1/40 | 7.9332e-07 | 2.3593 |
| 1/80 | 1.5510e-07 | 2.3547 |
| 1/160 | 3.2244e-08 | 2.2661 |

**Example 5:** Consider GFDEs [28] defined in Eq. (1) as,

$$^{*}\partial_t^{\alpha} u(x\,,t) - \frac{\partial^2 u(x\,,t)}{\partial^2 x} = \exp(x) t^4 \left(\frac{\Gamma(5+\alpha)}{24} - t^{\alpha}\right),\ x \in [0\,,1],\ t > 0, \tag{35}$$

with the initial and boundary conditions as $u(x\,,0) = 0,\ u(0\,,t) = t^{4+\alpha},\ u(1\,,t) = \exp(1) t^{4+\alpha}$.

The exact solution of the Eq. (1) under above condition is given by $u(x\,,t) = \exp(x) t^{4+\alpha}$. We take $z(t) = t,\ w(t) = 1$ and $0 \leq t \leq 1$. We solved Example 5 using the proposed method for $\alpha =\ 0.5$ and step sizes $\Delta x = 1/20000, \Delta t = 0.1$. The numerical results are shown in the Figures 5(a) and 5(b). The Figure 5(a) shows the exact solution while Figure 5(b) depicts numerical solution. We also solved

Example 5 with 1) varying both $\Delta x$ and $\Delta t$ with $\alpha = 0.15$; 2) varying $\Delta x$ and fixed $\Delta t = 1/2000$ with $\alpha = 0.4$ and 3) varying $\Delta t$ and fixed $\Delta x = 1/20000$ with $\alpha = 0.5$ which are shown through Table 5(a), 5(b) and 5(c) respectively. Tables 5(a), 5(b) and 5(c) establish that the proposed scheme is of second order of convergence. Table 5(b) and Table 5(c) show that the numerical scheme having second order of convergence in spatial direction and $3 - \alpha$ order of convergence in temporal direction respectively. Table 5(a), 5(b) and 5(c) are presented for the MAE and CO with $\alpha = 0.15, 0.4$ and $0.5$ respectively. It is clear from the Table 5(a), 5(b) and 5(c) our proposed numerical scheme is of high order of convergence and achieve the good accuracy. This example also concludes that the theoretical results with the proposed numerical results.

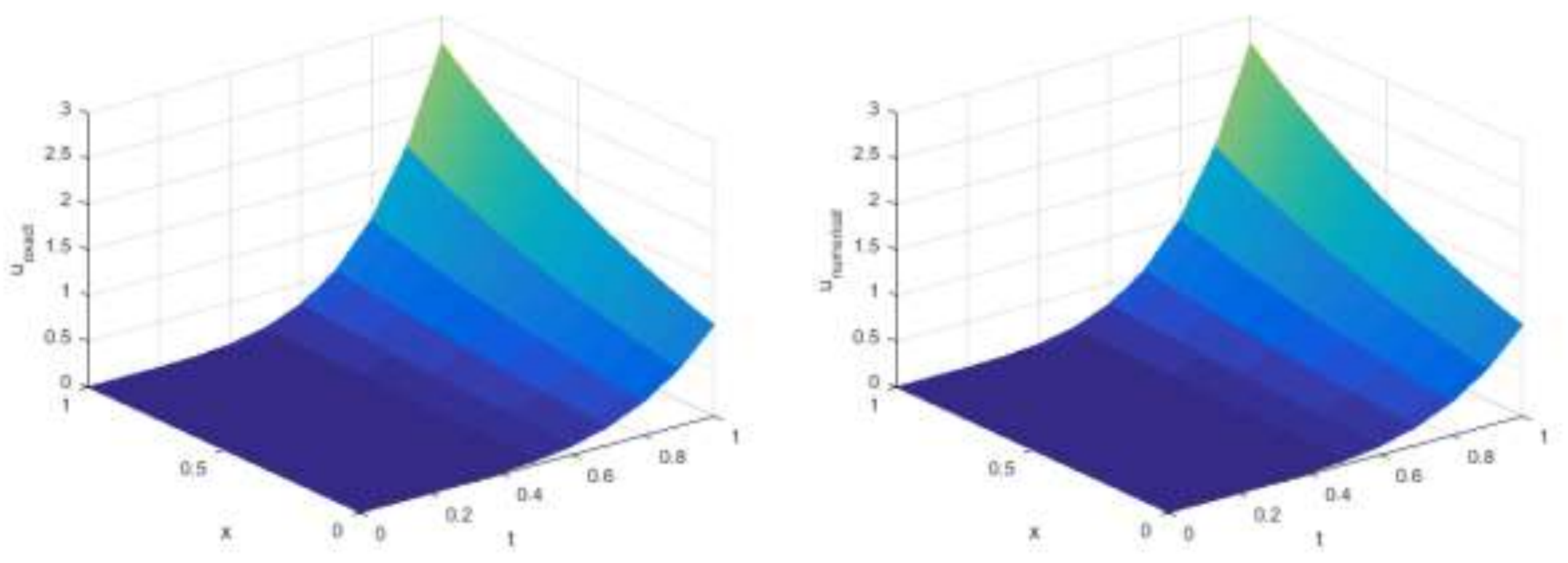

(a) Exact Solution. (b) $\Delta x = 1/20000, \Delta t = 0.1, \alpha = 0.5$.

**Fig 5:** Comparison of numerical and exact solutions of Example 5.

Table 5(a): MAE and CO for Example 5 with $\alpha = 0.15,\ \delta = 1\ z(t) = t,\ w(t) = 1$.

| $\Delta t$ | $\Delta x$ | MAE | CO |
|---|---|---|---|
| 1/10 | 1/10 | 3.6564e-04 | - |
| 1/20 | 1/20 | 7.2845e-05 | 2.3275 |
| 1/40 | 1/40 | 1.4938e-05 | 2.2858 |
| 1/80 | 1/80 | 3.2157e-06 | 2.2158 |
| 1/160 | 1/160 | 7.2456e-07 | 2.1499 |
| 1/320 | 1/320 | 1.6926e-07 | 2.0979 |

Table 5(b): MAE and CO for Example 5 with $\alpha = 0.4,\ \delta = 1\ z(t) = t,\ w(t) = 1\ \Delta t = 1/2000$.

| $\Delta x$ | MAE | CO |
|---|---|---|
| 1/10 | 1.4634e-04 | - |
| 1/20 | 3.6972e-05 | 1.9848 |
| 1/40 | 9.2475e-06 | 1.9993 |
| 1/80 | 2.3135e-06 | 1.9990 |
| 1/160 | 5.7981e-07 | 1.9964 |
| 1/320 | 1.4638e-07 | 1.9859 |

Table 5(c): MAE and CO for Example 5 with $\alpha = 0.5, \Delta x = \frac{1}{20000}, \delta = 1\ z(t) = t,\ w(t) = 1$ for $T = 1$.

| $\Delta t$ | MAE [29] | CO |
|---|---|---|
| 1/10 | 2.493639e-03 | - |
| 1/20 | 4.836648-04 | 2.37 |
| 1/40 | 9.015368e-05 | 2.42 |
| 1/80 | 1.644426e-05 | 2.45 |
| 1/160 | 2.972672e-06 | 2.47 |

## 6. Conclusion

We presented a high order scheme for solving the GFDEs defined in term of scale function and weight function in Caputo sense. We also proved the stability via convergence analysis of the proposed finite difference scheme. To check the performance of numerical scheme, we took five examples for performing numerical simulations. Simulation results show that numerical scheme is of high order and stable. We derived theoretically the stability of numerical scheme using the Fourier series method. Numerical results given in the form of Tables and Figures validate the presented numerical method. It was also noticed that the scheme works well with nonlinear scale and weight functions. The order of numerical scheme was investigated in temporal direction as $(\Delta t)^{3-\alpha}$ while in the spatial direction is $(\Delta x)^2$.